\documentclass{article}%
\usepackage{amssymb}
\usepackage{amsmath}
\usepackage{amsfonts}
\usepackage{graphicx}%
\providecommand{\U}[1]{\protect\rule{.1in}{.1in}}
\newtheorem{theorem}{Theorem}

\newtheorem{example}[theorem]{Example}

\newtheorem{remark}[theorem]{Remark}

\newenvironment{proof}[1][Proof]{\noindent\textbf{#1.} }{\ \rule{0.5em}{0.5em}}
\begin{document}

\title{ }

\begin{center}
\bigskip{\large Testing the fractional integration parameter revisited: a
Fractional Dickey-Fuller Test\bigskip}

Ahmed BENSALMA\footnote{High national school of statistic and applied economic
(ENSSEA), Algiers, Algeria, \ Email: \texttt{bensalma.ahmed@gmail.com }} and
Mohamed BENTARZI\footnote{ Faculty of Mathematics, University of Science and
Technology Houari Boumediene, Algiers, Algeria. , \ Email:
\texttt{bentarzimohamed@yahoo.fr }}\bigskip
\end{center}

\begin{description}
\item[Abstract] In this paper, in the first step, we show that the fractional
Dickey-Fuller test proposed by Dolado et al $\left[  10\right]  $ is useless
in practice. In the second step, we propose a new testing procedure for the
degree of fractional integration of a time series inspired on the unit root
test of Dickey-Fuller $\left[  7\right]  $. The composite null hypothesis is
that of $d\geq d_{0}$ against $d<d_{0}$. The test statistics is the same as in
Dickey-Fuller test using as output $\left(  1-L\right)  ^{d_{0}}y_{t}$ instead
of $\left(  1-L\right)  y_{t}$ and as input $\left(  1-L\right)  ^{-1+d_{0}%
}y_{t-1}$ and eventually some lag of $(1-L)^{d_{0}}y_{t}$ instead some lag of
$(1-L)y_{t}$, exploiting the fact that if $y_{t}$ is $I(d)$ then
$\Delta^{-1+d_{0}}y_{t}$ is $I(1)$ under the null $d=d_{0}$. If $d\geq d_{0}$,
using the generalization of Sowell's result $\left[  23\right]  $, we propose
a test based on the least favorable case, $d=d_{0}$, to control type $I$
error, and when $d<d_{0}$ we show that the tests statistics diverges to
$-\infty$, providing consistency. Through a simulation study, we show the good
performance of the test in terms of size and power. Finally, in order to show
how to use the new testing procedure, the test is applied to the well-known
Nelson and Plosser data.
\end{description}



\section{Introduction}

As the most popular long memory model and a useful extension of the classical
$ARIMA$ models, the fractionally integrated autoregressive moving average
($ARFIMA$) process, introduced by Granger and Jojeux $[11]$ and Hosking
$[13]$, has seen a considerable interest in the past three decades and has
been widely applied in many fields like hydrology, economics and finance. The
$ARFIMA$ process generalizes the standard linear $ARIMA(p,d,q)$ model by
permitting to the degree of integration $d$ to be non-integer. Compared with
the standard $ARMA$ and $ARIMA$ specifications, the $ARFIMA$ generalization
provides a more flexible framework in modelling the long range dependence,
where a special role is played by the fractional differencing parameter $d$
whose precise determination is very important in applied work.

The stationary and invertible Fractional $ARIMA(p,d,q)$ processes, is defined
as the following%
\[
\left(  1-\sum_{i=1}^{p}\phi_{i}L^{i}\right)  \left(  1-L\right)  ^{d}%
y_{t}=\left(  1-\sum_{i=1}^{q}\theta_{i}L^{i}\right)  u_{t}\text{, \ \ }t\in%
\mathbb{Z}
\text{, \ \ }-0.5<d<0.5,
\]
where $L$ is the backshift operator and $u_{t}$ are independently and
identically distributed ($i.i.d$) random variables with zero mean and finite
variance; $\left(  1-\sum_{i=1}^{p}\phi_{i}L^{i}\right)  $ and $\left(
1-\sum_{i=1}^{q}\theta_{i}L^{i}\right)  $ are polynomial functions of $L$ with
order $p$ and $q$, and both of them have only roots outside the unit circle.
The fractional difference operator $(1-L)^{d}$ is defined by its Maclaurin
series (by its binomial expansion, if $d$ is an integer):%
\begin{align*}
(1-L)^{d} &  =\sum_{j=0}^{\infty}\frac{\Gamma(-d+j)}{\Gamma(-d)\Gamma
(j+1)}L^{j}\\
&  =\sum_{j=0}^{\infty}\frac{(-d)(-d+1)\cdots(-d+j-1)}{j!}L^{j},
\end{align*}
where%
\[
\Gamma(z)=\left\{
\begin{array}
[c]{cc}%
\int_{0}^{\infty}s^{z-1}e^{-z}ds & \text{If }z>0\\
\infty & z=0.
\end{array}
\right.
\]
If $z<0$, $\Gamma(z)$ is defined by the recursion formula $z\Gamma
(z)=\Gamma(z+1)$.

In recent years, an increasing effort has been made to establish reliable
testing procedures to determine whether or not an observed time series is
fractionally integrated. In particular, there has been a considerable interest
in generalizing the familiar Dickey-Fuller test by taking into account the
fractional integration order. It is well documented that the power of
Dickey-Fuller [$DF$] type tests against alternatives of fractional integration
is low (see Sowell $\left[  23\right]  $; Diebold and Rudebusch $\left[
9\right]  $; Hassler and Wolters $\left[  12\right]  ;$ Kr\"{a}mer $\left[
14\right]  $). This motivated the development of powerful tests against
fractional alternatives. Robinson $\left[  20\right]  $ pioneered an
integration test constructed from the Lagrange Multiplier [$LM$] principle,
which was proven by Robinson $\left[  21\right]  $ to be locally the most
powerful under Gaussianity. The test has been further studied and modified by
Agiakloglou and Newbold $\left[  1\right]  $, Tanaka $\left[  24\right]  $.
Tanaka $\left[  24\right]  $ showed, through simulation experiments, that the
$LM$ tests have serious size distortion. Another serious criticism addressed
to the $LM$ tests is that, by working under the null hypothesis, it does not
yield any direct information about the correct long-memory parameter $d$, when
the null is rejected (Candelon, Gil Alana $\left[  6\right]  $).

More recently, Dolado et \textit{al} $\left[  10\right]  $ introduced a
fractional integration test (henceforth $DGM$ test) based on an auxiliary
regression for the null of unit root $(H_{0}:d=1)$ against the alternative of
fractional integration ($H_{1}:d=d_{1}$, $d_{1}<1$). Their proposed test
reduces to the standard Dickey-Fuller test when $d_{1}=0$ while under the null
and when $d_{1}$ known, the statistic in the corresponding regression model
depends on a fractional Brownian motion if $0\leq d_{1}<0.5$.\ Further, the
$DGM$ test was refined by Lobato and Velasco $(\left[  16\right]  ,\left[
17\right]  )$ using the same null and alternative hypotheses.

While the $DGM$ test represents a useful generalization of the Dickey-Fuller
test in the presence of a fractionally integrated alternative, it might give
arbitrary conclusions when the true $d$ is not present neither in the null nor
in the alternative, because the auxiliary regression model, they used, depends
on the null and alternative (i.e. $1$ and $d_{1}$)\emph{.} Indeed, through
some simulation experiments we conduct, it may be seen (see Table $1$ below)
that the $DGM$ test performs somewhat badly in the case where the parameter
$d$ is wrongly specified under the null and alternative. In such situation
three cases can arise: the case where the null is true, the case where the
alternative is true and the case where neither the null nor the alternative is true.

In this paper, we propose an alternative test for the fractional parameter
$d,$ inspired by the unit root test of Dickey-Fuller $\left[  7\right]  $. The
composite null hypothesis is that of $d\geq d_{0}$ against $d<d_{0}$. The test
statistics is the same as in Dickey-Fuller test using as input $\Delta
^{-1+d_{0}}y_{t-1}$ instead of $y_{t-1}$, exploiting the fact that if $y_{t}$
is $I(d)$ then $\Delta^{-1+d_{0}}y_{t}$ is $I(1)$ under the null $d=d_{0}$. If
$d\geq d_{0}$, using the generalization of Sowell's result $\left[  23\right]
$, we propose a test based on the least favorable case, $d=d_{0}$, to control
type $I$ error, and when $d<d_{0}$ we show that the tests statistics diverges
to $-\infty$, providing consistency. Clearly such testing procedure is
conceptually attractive since, first, the hypotheses we consider are rather
composite-versus-composite $(H_{0}:d\geq d_{0}$ against $H_{1}:d<d_{0}$)
resulting in a dichotomic choice which excludes the third case. Second, by the
choice of a suitable regression model, $\Delta^{d_{0}}y_{t}=\rho
\Delta^{-1+d_{0}}y_{t-1}+\varepsilon_{t}$, the usual statistics $t_{\widehat
{\rho}_{n}}$ or $n\widehat{\rho}_{n}$ have the same asymptotic distribution as
the Dickey-Fuller test.\ This is because the maximum probability of rejecting
the null hypothesis i.e. $\alpha=Sup_{d\geq d_{0}}P(reject$ $H_{0})$, level of
the test, is reached when $d=d_{0}$. So, the standard Dickey-Fuller table may
be used for our test without an extra-effort i.e. without using the tabulated
values of a fractional Brownian motion.

Before going through the topic, it is important to precise certain essential
points, which may facilitate the reading of this paper. The main theme of our
article is how to extend the familiar Dickey-Fuller $\left[  7\right]  $ type
tests for unit root ($I(1)$ against $I(0)$) by embedding the case $d=0$ and
$d=1$ in continuum of memory properties (i.e. $d\in%
\mathbb{R}
$). Such extension has already been discussed by Dolado et Al $[10]$. In our
paper, we show, in the first step, that the $DGM$ approach is not the best and
adequate way to extend the Dickey-Fuller test by taking into account the
fractional case. In the second step, we provide how to extend adequately the
standard Dickey-Fuller test $\left[  7\right]  $ by taking into account the
fractional case.

In order to expose clearly the alternative approach and permit the careful
comparison with the $DGM$ approach, we choose to use a simple framework like
$ARFIMA(0,d,0)$($\equiv FI(d)$) process. The case, where the errors are
autocorrelated, deserves that one devotes another paper, by taking into
account the seminal work of Said and Dickey $\left[  22\right]  $ and Phillips
$\left[  19\right]  $. Our approach is based on the following forth points:

\begin{description}
\item[1] Using the composite hypothesis $H_{0}:d\geq d_{0}.$

\item[2] If $y_{t}\leadsto I(d_{0})$ then $(1-L)^{-1+d_{0}}y_{t}\leadsto
I(1).$

\item[3] Testing the composite null hypothesis is based upon testing the
statistical significance of the coefficient $\phi$ (or $\rho=\phi-1$) in the
regression model $\Delta^{-1+d_{0}}y_{t}=\phi\Delta^{-1+d_{0}}y_{t-1}%
+\varepsilon_{t}.$

\item[4] The level of the test $\alpha=Sup_{d\geq d_{0}}P(reject$
$H_{0})=P(reject$ $H_{0}|d=d_{0}).$
\end{description}

In order to highlight these four important points and not to overlook them
into a mid-general framework, the case where the errors are correlated, will
not be pursued in this paper. However, I provide (see Appendix $2$) some
discussions when the process $\left\{  y_{t}\text{, \ }t\in%
\mathbb{Z}
\right\}  $ is generated by%
\[
y_{t}=\mu(t)+FI(d),
\]
with $\mu(t)$ being a vector of deterministic functions like a constant or
time trend.

Another reason that led us to choose a simple theoretical framework is to
highlight the importance of considering correctly, some basic rules of the
testing statistical hypothesis theory. In this paper, we focus on the
importance to consider the statistic of the test, exclusively, deduced under
the null hypothesis (see section $2$, for more details).

The rest of this paper is organized as follows. In section $2$, to highlight
the contribution of our approach, we first give some comments on the $DGM$
approach. Then in section $3$, we define in a simple framework our test and in
particular the auxiliary regression model used to test the null. Moreover, the
main results on the asymptotic distribution under the null and alternative
composite hypothesis are given. Section $4$ explores a theoretical study about
the size and power of our proposed $F$-$DF$ test. Furthermore, Monte-Carlo
simulation experiments are undertaken in order to support the analytical
results and in particular to confirm that the proposed test is robust to any
misspecification of the order of integration parameter $d$. In Section $5$, we
present empirical applications by revisiting Nelson-Plosser data. It is
important to note that the empirical application is made only to explain how
to use the new testing procedure (The reader should not understand this
application as to provide a new evidence for the order of integration of the
Nelson and Plosser data). Because, as it has been mentioned previously, the
data generating process adopted in this paper is restrictive). Finally, the
proofs of the main results presented in Section $2$ are left to the appendix
$1$ and some discussions when the process $\left\{  y_{t}\text{, \ }t\in%
\mathbb{Z}
\right\}  $ is generated by some deterministic trend plus $FI(d)$ process are
left to the Appendix $2$.

\section{Fractional Dickey-Fuller testing: the $DGM$ approach}

\subsection{Hypotheses and the auxiliary regression model}

Dolado, Gonzalo, and Mayoral [$DGM$] $\left[  10\right]  $ introduced a test
based on an auxiliary regression for the null of unit root against the
alternative of fractional integration. The fractional Dickey-Fuller ($F$-$DF$)
test considered by $DGM$ $\left[  10\right]  $, in the basic framework, is
described as follows.

Let $\left\{  y_{t}\right\}  _{t=1}^{n}$ a series generated from the
fractionally integrated process ( $FI(d)$ in short) given by%
\begin{equation}
(1-L)^{d}y_{t}=u_{t},\text{ }t\in%
\mathbb{Z}
, \tag{$2.1$}%
\end{equation}
where $d\in%
\mathbb{R}
$ is the true order of integration and, $\left\{  u_{t},t\in%
\mathbb{Z}
\right\}  $ is an $iid$ innovation with mean zero and variance $\sigma_{u}^{2}
$. For the data generating process ($DGP$) $(2.1)$, $DGM$ $\left[  10\right]
$ propose to test the following hypotheses,%
\begin{equation}
H_{0}:d=d_{0}\text{ against }H_{1}:d=d_{1},\text{with }d_{1}<d_{0},
\tag{$2.2$}%
\end{equation}
by means of the $t$ statistic of the coefficient of $\Delta^{d_{1}}y_{t-1}$ in
the ordinary least squares ($OLS$) regression%
\begin{equation}
\Delta^{d_{0}}y_{t}=\rho\Delta^{d_{1}}y_{t-1}+\varepsilon_{t,}\text{
(}t=1,\cdots,n\text{),} \tag{$2.3$}%
\end{equation}
where $\Delta=1-L$. $DGM$ explains the choice of the auxiliary regression
model $(2.3)$ by arguing that, in the simple Dickey-Fuller test, to test the
hypotheses
\begin{equation}
H_{0}:d=1\text{ against }H_{1}:d=0, \tag{$2.4$}%
\end{equation}
the maintained regression model is:%
\begin{equation}
\Delta y_{t}=\rho y_{t-1}+\varepsilon_{t,}\text{ \ }t=1,\cdots,n, \tag{$2.5$}%
\end{equation}
where $\varepsilon_{t}\sim iid(0,\sigma_{\varepsilon}^{2})$. If $y_{t}$ is
$I(1)$, then the regression $(2.5)$ is unbalanced in the sense that the orders
of integration of the regressand and the regressor are different, being $I(0)$
and $I(1)$ respectively. After this, $DGM$ claim that in the simple
Dickey-Fuller test the null hypothesis $H_{0}:d=1$ correspond to the regressor
$\Delta^{1}y_{t}$ and the alternative $H_{1}:d=0$ correspond to $\Delta
^{0}y_{t-1}=y_{t-1}$. This leads them to consider that the null hypothesis
$H_{0}:d=d_{0}$ correspond to $\Delta^{d_{0}}y_{t}$ and the alternative
$H_{1}:d=d_{1}$ correspond to $\Delta^{d_{1}}y_{t-1}$. In the following we
show that the interpretation, made by $DGM$, in the use of the model $(2.5)$,
in the simple Dickey-Fuller test is incorrect. In fact, the standard Dickey
and Fuller test is not based directly on the regression model $(2.5)$. The
hypotheses $(2.4)$ are based on the following regression model%
\begin{equation}
y_{t}=\phi y_{t-1}+\varepsilon_{t}\text{ \ }t=1,\cdots,n, \tag{$2.6$}%
\end{equation}
which is equivalent to the regression model $(2.5)$, with $\rho=\phi-1$. The
regression model $(2.6)$ is balanced in the sense that the regressand and the
regressor have the same order of integration which is equal $1$ under the
null. The scheme $1$ and $2$ summarize, respectively, the incorrect and
correct interpretation in the use of the model $(2.5)$, in the simple
Dickey-Fuller test.\bigskip

\begin{center}%
\begin{tabular}
[c]{|c|c|}\hline
\multicolumn{2}{|c|}{\textbf{Scheme }$\mathbf{1}$: Incorrect interpretation,
in the use of model $(2.5)$,}\\
\multicolumn{2}{|c|}{in the simple Dickey-Fuller test}\\\hline
$DGM$ interpretation, in the use of the model & The use of the incorrect\\
$(2.5)$, in the simple Dickey-Fuller test & interpretation in fractional
case\\\hline
${\LARGE \Delta}^{\overset{\overset{{\LARGE H}_{0}{\LARGE :d=1}}{\downarrow}%
}{1}}{\LARGE y}_{t}{\LARGE =\rho\Delta}^{\overset{\overset{{\LARGE H}%
_{1}{\LARGE :d=0}}{\downarrow}}{0}}{\LARGE y}_{t-1}{\LARGE +\varepsilon}_{t}$
& ${\LARGE \Delta}^{\overset{\overset{{\LARGE H}_{0}{\LARGE :d=d}_{0}%
}{\downarrow}}{{\LARGE d}_{0}}}{\LARGE y}_{t}{\LARGE =\rho\Delta}%
^{\overset{\overset{{\LARGE H}_{1}{\LARGE :d=d}_{1}}{\downarrow}}%
{{\LARGE d}_{1}}}{\LARGE y}_{t-1}{\LARGE +\varepsilon}_{t}$\\\hline
\end{tabular}

\bigskip%

\begin{tabular}
[c]{|c|c|}\hline
\multicolumn{2}{|c|}{\textbf{Scheme }$\mathbf{2}$: Correct interpretation, in
the use of the models $(2.5)$ and $(2.6)$,}\\
\multicolumn{2}{|c|}{in the simple Dickey-Fuller test}\\\hline
Correct interpretation, in use & Correct interpretation, in the use\\
of the model $(2.5)$ & of the model $(2.6)$\\\hline
$\overset{\overset{{\LARGE I(0)}}{{\LARGE under}\text{ }{\LARGE H}_{0}}%
}{\overbrace{{\LARGE \Delta y}_{t}}}{\LARGE =\rho}\overset{\overset
{{\LARGE I(1)}}{{\LARGE under}\text{ }{\LARGE H}_{0}}}{\overbrace
{{\LARGE y}_{t-1}}}{\LARGE +\varepsilon}_{t}$ & $\overset{\overset
{{\LARGE I(1)}}{{\LARGE under}\text{ }{\LARGE H}_{0}}}{\overbrace
{{\LARGE y}_{t}}}={\LARGE \phi}\overset{\overset{{\LARGE I(1)}}{{\LARGE under}%
\text{ }{\LARGE H}_{0}}}{\overbrace{{\LARGE y}_{t-1}}}+{\LARGE \varepsilon
}_{t}$\\\hline
\end{tabular}

\bigskip
\end{center}

The scheme $2$, indicates that the choice of the model $(2.5)$ or equivalently
$(2.6)$ is based only on the null hypothesis. The scheme $1$, indicates that
the incorrect interpretation leads $DGM$ to consider a regression model based
on the null and alternative in the fractional case. In fact, for the
regression model $(2.3)$ we have, under the null ($H_{0}:d=d_{0}$)%
\[
\Delta^{d_{0}}y_{t}\leadsto I(0)\text{ \ \ and \ \ \ }\Delta^{d_{1}}%
y_{t}\leadsto I(d_{0}-d_{1}),
\]
and under the alternative ($H_{0}:d=d_{1}$),%
\[
\Delta^{d_{0}}y_{t}\leadsto I(d_{1}-d_{0})\text{ \ \ and \ \ \ }\Delta^{d_{1}%
}y_{t}\leadsto I(0).
\]

As a result, the authors are locked into the trap set by this semblance of
analogy. There are other inconsistencies in the use of statistical concepts.
Someone, can easily feel, throughout the reading of the article $[10]$, the
efforts granted by the authors to justify inconsistencies. It would be long to
enumerate all the inconsistencies in the $DGM$ test procedure.

\subsubsection{Unit root test against fractional alternatives and its
asymptotic Properties}

To study the performances of their procedure in terms of power and size, $DGM
$ $\left[  10\right]  $ consider only the particular case,%
\begin{equation}
H_{0}:d=1\text{ against }H_{1}:d=d_{1}, \tag{$2.7$}%
\end{equation}
by means of the $t$-statistic of the coefficient $\Delta^{d_{1}}y_{t-1}$, in
the ordinary least squares ($OLS$) regression%
\begin{equation}
\Delta^{1}y_{t}=\rho\Delta^{d_{1}}y_{t-1}+\varepsilon_{t}. \tag{$2.8$}%
\end{equation}
The $t$-ratio, $t_{\widehat{\rho}}(d_{1})$, is given by%
\[
t_{\widehat{\rho}}(d_{1})=\frac{\sqrt{n}\sum_{t=2}^{n}\Delta y_{t}%
\Delta^{d_{1}}y_{t-1}}{\sqrt{\sum_{t=2}^{n}\left(  \Delta y_{t}-\widehat{\phi
}\Delta^{d_{1}}y_{t-1}\right)  ^{2}\sum_{t=2}^{n}\left(  \Delta^{d_{1}}%
y_{t-1}\right)  ^{2}}}.
\]

\begin{theorem}
(DGM [10]. Under the null hypothesis that $y_{t}$ is a random walk, the
asymptotic distribution of $t_{\widehat{\rho}}(d_{1})$ is given by%
\[
t_{\widehat{\rho}}(d_{1})\overset{L}{\longrightarrow}\frac{\int w_{-d_{1}%
}(r)dB(r)}{\left(  \int w_{-d_{1}}^{2}(r)dr\right)  ^{1/2}}\text{ if }0\leq
d_{1}<0.5,
\]
and%
\[
t_{\widehat{\rho}}(d_{1})\overset{L}{\longrightarrow}N(0,1)\text{ if }0.5\leq
d_{1}<1.
\]
where $w_{-d_{1}}(\cdot)$ is fractional Brownian motion.
\end{theorem}

\begin{proof}
See DGM [10]
\end{proof}

Theorem $1$, shows that under the null the asymptotic distribution of
$t$-statistic depends on fractional Brownian motion if $0\leq d_{1}<0.5$ and
$t_{\widehat{\rho}}(d_{1})\longrightarrow N(0,1)$ if $0.5\leq d_{1}<1$. These
asymptotic distributions are different from those derived by Dickey and Fuller
$\left[  7\right]  $ which depend only on standard Brownian motion. The
implementation of $DGM$ $\left[  10\right]  $ test would require tabulation of
the percentiles of the functional of Brownian motion, which imply that
inference on the presence of unit root would be conditional on $d_{1}$. But
given the well-known difficulties in estimating the order of fractional
integration in finites samples, thus the test might suffer from
misspecification (i.e. the parameter $d$ is wrongly specified)

Under $H_{0}:d=1$, we have $Cov(\Delta^{d_{0}}y_{t},\Delta^{d_{1}}y_{t-1})=0$
and under the alternative we have $Cov(\Delta^{d_{0}}y_{t},\Delta^{d_{1}%
}y_{t-1})=\sigma_{u}^{2}(-1+d_{1})<0$. Thus $DGM$ build the decision rule as
follows,%
\begin{equation}
\left\{
\begin{array}
[c]{c}%
\begin{array}
[c]{cc}%
H_{0}:d=1\text{ is accepted} & \text{if }\rho=0
\end{array}
\\%
\begin{array}
[c]{cc}%
H_{0}:d=d_{1}\text{ is accepted} & \text{if }\rho<0
\end{array}
\end{array}
\right.  \tag{$2.9$}%
\end{equation}

The hypotheses $(2.7)$ based on the regression model $(2.8)$ and the decision
rule $(2.9)$ is called by their authors "Fractional Dickey and Fuller Test".

\begin{remark}
Why does DGM consider only the case $H_{0}:d=1$? To respond this question, let
us consider the case $H_{0}:d=0.5$. Since $\Omega(d)=[0,1]$, the set of
alternatives values of $d$ is%
\[
\Omega_{1}(d)=\left[  0,1/2\right[  \cup\left]  1/2,1\right]  .
\]
In this example, we have two cases. The first case, is given by $d_{0}=1/2$
and $d_{1}\in\left[  0,1/2\right[  $ and the decision rule is based on
"($\rho=0$ or $\rho<0$). The second case is given by $d_{0}=1/2$ and $d_{1}%
\in\left]  1/2,1\right]  $. In this case it is easy to show that the decision
rule is based on "($\rho=0$ or $\rho>0$)", because we have%
\begin{align*}
\text{under the null }Cov(\Delta^{0.5}y_{t},\Delta^{d_{1}}y_{t-1})  &  =0,\\
\text{under the alternative }Cov(\Delta^{0.5}y_{t},\Delta^{d_{1}}y_{t-1})  &
=\sigma_{u}^{2}(-0.5+d_{1})>0.
\end{align*}

\end{remark}

\begin{remark}
We can suggest another D.G.M. type test. Indeed, since DGM would only use the
decision rule $(2.9)$, they would have been better advised if they had thought
about testing hypothesis%
\begin{equation}
H_{0}:d=d_{0}\text{ against }H_{1}:d=0,\text{ with }d_{0}>0. \tag{$2.10$}%
\end{equation}
This choice can be justified by the integration order $d>0$ of the majority of
economic series. By using the scheme 1, we can deduce that to test $(2.10)$ we
must use the regression model%
\begin{equation}
\Delta^{d_{0}}y_{t}=\rho y_{t-1}+\varepsilon_{t}. \tag{$2.11$}%
\end{equation}
With the test hypotheses (2.10) and the regression model (2.11) we can use the
decision rule "($\rho=0$ ou $\rho<0$), since the case $d<0$ is excluded.
\end{remark}

\subsubsection{Power and size of DGM's FDF test.}

The problem with the test based on the hypotheses $(2.7)$ and regression model
$(2.8)$ and the test suggested above, based on $(2.10)$ and $(2.11)$ are
useless in practice. The problem with the $DGM$ type tests is that they are
based on a choice of two possible orders of integration $d_{0}$ and $d_{1}$,
of which the true order can be different either in the null or in the
alternative. In fact, in the fractional integration case, there is a continuum
of possible orders of integration. This would make the simple-versus-simple
hypothesis invalid, particularly if the auxiliary regression model, used for
the test, is based on the null and alternative. For instance, in the $DGM$
test one of the following three cases holds:

\begin{itemize}
\item $d=d_{0}$,

\item $d=d_{1}$,

\item $d\neq d_{0}$ and $d\neq d_{1}$.
\end{itemize}

The third case causes serious troubles in practice, particularly, if the
statistic of the test depends on null and alternative hypothesis. When
$d_{0}=1$, in the first two cases, Dolado et al $\left[  10\right]  $ showed
by means of a simulation study that their test procedure has a good
performance in terms of power and level. For the third case, Dolado et al
$\left[  10\right]  $ studied the effect of hypotheses misspecification by
considering the deviations from the true value $d_{1}$ with size $\pm0.1$,
$\pm0.2$ and $\pm0.3$. In the following; however, we replicate the simulation
results of Dolado et al $\left[  10\right]  $ and present them more clearly by
using a single table. We generate $1000$ series from the data generating
process $(2.1)$ with sample size $n=100$. The first column of Table $1$ gives
the true values of the parameter $d$ while the second line shows the values of
$d_{1}$ specified under the alternative. The first line gives the tabulated
values by $DGM$ (see Dolado et al $\left[  10\right]  $, table $X$\ page
$2003$). The last line of Table $1$ represents the performance of the $DGM$
test in terms of level, i.e. the percentage of rejection of the null, when it
is true ($\alpha$), while the main diagonal represents the performance of the
$DGM$ test in terms of power i.e. the percentage of acceptance of the
alternative hypothesis when it is true, ($1-\beta$). $\alpha$ and $\beta$ are
respectively the type $I$ and the type $II$ errors, defined by%
\[
\alpha=P(\text{reject }H_{0}|d=1\text{) \ \ \ and \ \ \ }\beta=P(\text{reject
}H_{1}|d=d_{1}\text{).}
\]
The other values in the table are the percentage of acceptance of the
alternative hypothesis when both the null and alternative are false i.e. when
the value of $d$ is wrongly specified. In fact, these values represents
another type of errors, namely%
\[
P_{d\neq d_{1}}\left(  \text{Accept }H_{1}|d\neq1\text{ and }d\neq
d_{1}\right)  .
\]
When performing a test one may arrive at the correct decision, or one may
commit one of two errors: rejecting the null hypothesis when it is true (type
$I$ error, or error of the first kind) or accepting it when it is false (type
$II$ error or error of the second kind). In statistical testing theory, there
is no place for type $III$ error (or error of the third kind). This anomaly is
the consequence of the choice of inappropriate auxiliary regression model,
which depends on the null and alternative. From Table $1$, it may be easily
observed that when the true $d$ is well specified, the $DGM$ test has a good
performance in terms of power and level. However, in the case where the true
value of $d\in\left[  0,1\right]  -\left\{  1,d_{1}\right\}  $, the
conclusions of the test are somewhat arbitrary. For example, when $d=0.3$, the
percentage of acceptance of the alternative is equal $100\%$ regardless of the
alternative hypothesis. In other word, if the process $y_{t}$, is fractionally
integrated of order $d=0.3$ (i.e. stationary stationary process), the table
$1$, show that for $H_{0}:d=1$ against $H_{1}:d=0.7$, we have%
\[
P_{d=0.3}\left(  \text{Accept }H_{1}:d=0.7|d\neq1\text{ and }d\neq
d_{1}\right)  =1.
\]
This example shows clearly that the risk to specify the stationary process as
a nonstationary process is high.

\begin{figure}[ptbh]
\begin{center}
\includegraphics[width=12cm]{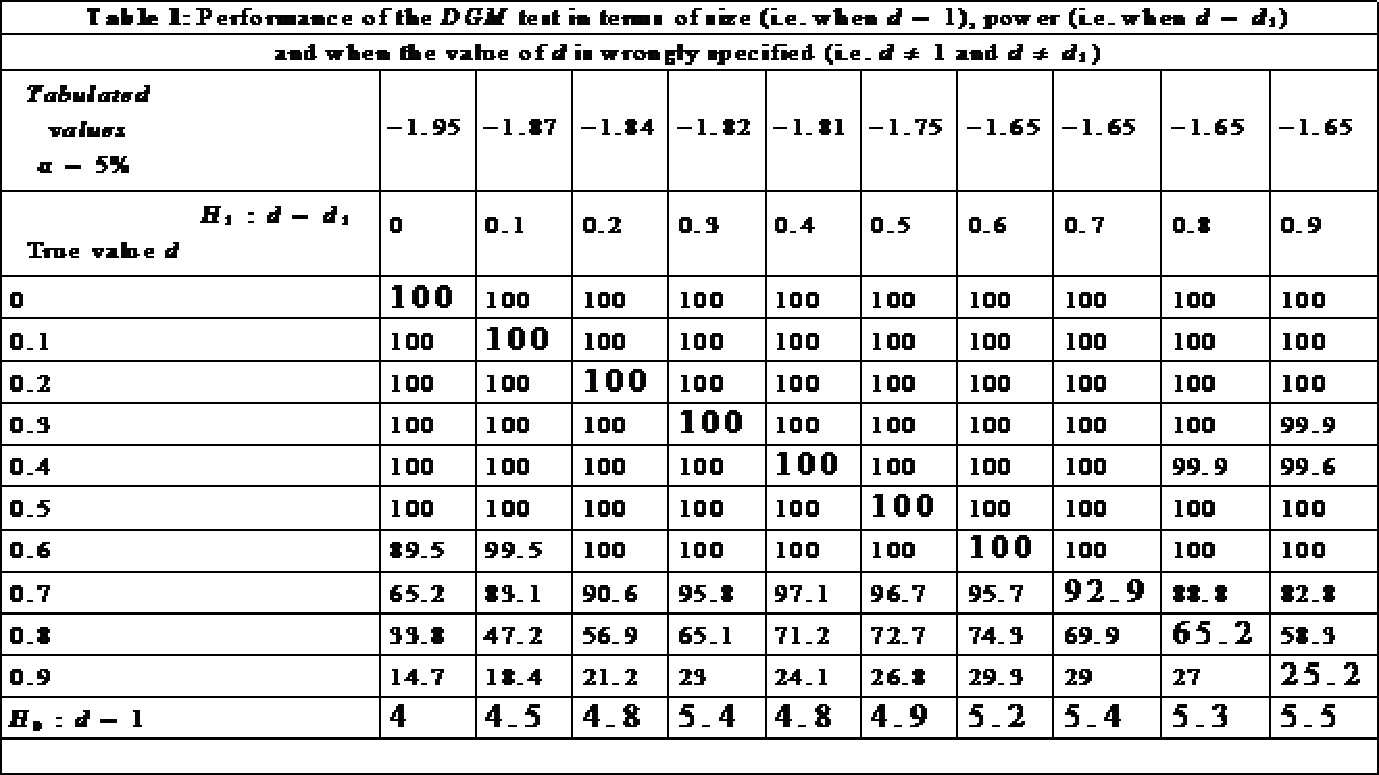}
\end{center}
\end{figure}$LV$ $\left[  17\right]  $ argue that $\Delta^{d_{1}}y_{t-1}$ is
not the best class of regression one can choose and propose another auxiliary
regression model for the test $(2.4)$. In the case $d_{0}=1,$ they propose to
test $(2.4)$ by using the following auxiliary model%
\[
\Delta y_{t}=\phi_{2}z_{t-1}(d_{1})+\varepsilon_{t,}\text{ (}t=1,\cdots
,n\text{),}
\]
where%
\[
z_{t-1}(d_{1})=\left(  \frac{\Delta^{d_{1}-1}-1}{1-d_{1}}\right)  \Delta
y_{t}.
\]
The same criticisms can be formulated concerning test concerning $LV$ test.
The $DGM$ $\left[  10\right]  $ and $LV$ $\left[  17\right]  $ tests present
an analogy with the original Dickey-Fuller test, but can not be considered as
a generalization of the familiar Dickey-Fuller test in the sense that the
conventional $I(1)$ vs $I(0)$ framework is recovered (for the $DGM$ test the
conventional framework is recovered only if $d_{0}=1$ and $d_{1}=0$%
).\emph{\ }The implementation of $DGM$ $\left[  10\right]  $ test would
require tabulations of the percentiles of the functional of fractional
Brownian motion, which imply that the inference on the presence of unit root
would be conditional on $d_{1}$, and thus might suffer from misspecification
resulting from errors in specifying the fractional parameter $d$. When $d_{1}
$ is not taken to be known a priory, a pre-estimation of it is needed to
implement the test. In this case, we can perform the test only if the
estimator of $d_{1}$ ($\widehat{d}_{1}$) is sufficiently close to unity (see
$DGM$ $\left[  10\right]  $ for more details). Indeed, the table $1$, show
that the $DGM$ test have "a realistic" behavior (i.e. likely results) in terms
of size and power only when the true value of $d$ is close to $1$ (see for
instance in table $1$:, $d=0.8$ and $d=0.9$). This is why, $DGM$ recommend to
use an estimator of $d_{1}$ that originates from the trimming rule%
\[
\widehat{d}_{1}=\left\{
\begin{array}
[c]{c}%
\begin{tabular}
[c]{cc}%
$\widehat{d}_{n}$, & if $\widehat{d}_{n}<1-c$%
\end{tabular}
\\%
\begin{tabular}
[c]{cc}%
$1-c$ & if $\widehat{d}_{n}\geq1-c$%
\end{tabular}
\end{array}
\right.  ,0<c<0.5,
\]
where $\widehat{d}_{n}$ is any $\sqrt{n}$-consistent estimator of $d_{1}$, for
example $DGM$ select $c=0.02$ in their simulation experiments. This rule
makes\emph{\ }this test more vague in how to use it in practice.

To extend adequately the standard Dickey-Fuller test $[7]$, we propose a new
test based on mutually exclusive and complementary null, alternative
hypotheses and a suitable auxiliary regression model.

\section{Fractional Dickey-Fuller testing: an alternative approach}

\subsection{Hypotheses, the auxiliary regression model and asymptotic under
the null and the alternative}

In this section, we deal with a series $\left\{  y_{t}\right\}  _{t=1}^{n}$
generated from the fractionally integrated model, $FI(d)$, given by $(2.1)$,
where the order $d$ is any real number. Under this setting, we propose to test
the following hypotheses\footnote{The special case of hypothesis testing
$H_{0}:d\geq1$ against $H_{1}:d<1$ was presented at ICMSAO'13 Conference,
Hammamet, Tunisia, 28--30 April 2013, in the paper entitled "A consistent test
for unit root against fractional alternative". Expanded version of this paper
forthcoming in Inderscience journal "International Journal of operational
research "}\footnote{This paper is an expanded version of the paper entitled
"New fractional Dickey-Fuller test" presented at ICMSAO'15 conference,
Istanbul, May 27-29,2015} :%
\begin{equation}
H_{0}:d\geq d_{0}\ \ \ \ against\ \ \ \ H_{1}:d<d_{0}. \tag{$3.1$}%
\end{equation}
Our proposal is based upon testing the statistical significance of the
coefficient $\phi$ (or $\rho=\phi-1$) in the following regression model,%
\begin{equation}
\Delta^{-1+d_{0}}y_{t}=\phi\Delta^{-1+d_{0}}y_{t-1}+\varepsilon_{t},
\tag{$3.2$}%
\end{equation}
or equivalently%
\begin{equation}
\Delta^{d_{0}}y_{t}=\rho\Delta^{-1+d_{0}}y_{t-1}+\varepsilon_{t}, \tag{$3.3$}%
\end{equation}
where $\rho=\phi-1$ and $\left\{  \varepsilon_{t}\right\}  _{t=1}^{n}$ are the
residuals. The most important idea behind the choice of the framework above is
that%
\[
\text{when }d=d_{0}\text{, \ \ }x_{t}=\Delta^{-1+d_{0}}y_{t}\text{, is
integrated of order }1\text{.}
\]
More generally,%
\[
x_{t}\text{ is integrated of order }1+d-d_{0},
\]
with%
\[
\left\{
\begin{array}
[c]{c}%
1+d-d_{0}\geq1\text{, if }d\geq d_{0},\\
1+d-d_{0}<1\text{, if }d<d_{0}.
\end{array}
\right.
\]
Before stating the main results of this paper, we give some technical tools
that we need in the sequel. Let $\eta_{t}=(1-L)^{-\delta}u_{t}$, with
$\delta\in\left]  -0.5,0.5\right]  $ and $u_{t}$ defined as above. Let
$\sigma_{S}^{2}=var(S_{n})$, where $S_{t}=\sum_{j=1}^{t}\eta_{j}$. When
$\left\vert \delta\right\vert <\frac{1}{2}$, we have (see Sowell $\left[
23\right]  )$%
\begin{equation}
\underset{n\rightarrow\infty}{\lim}n^{-1-2\delta}\sigma_{S}^{2}=\frac
{\sigma_{\varepsilon}^{2}\Gamma(1-2\delta)}{(1+2\delta)\Gamma(1+\delta
)\Gamma(1-\delta)}\equiv\kappa_{\eta}^{2}(\delta), \tag{$3.4$}%
\end{equation}
where $\Gamma(\cdot)$ denotes the Gamma or generalized factorial function. For
the case $\delta=\frac{1}{2}$, (see Liu, $\left[  15\right]  )$%
\begin{equation}
\underset{n\rightarrow\infty}{\lim}(n^{-2}\log^{-1}n)\sigma_{S}^{2}%
=\frac{2\sigma_{\varepsilon}^{2}}{\pi}\equiv\kappa_{\eta}^{2}(\frac{1}{2}).
\tag{$3.5$}%
\end{equation}
Furthermore, under the following additional assumption $E\left\vert
u_{t}\right\vert ^{a}<\infty,$ for some $a\geq\max\left\{  4,\frac{-8\delta
}{1+2\delta}\right\}  $, the following useful results apply:%
\begin{equation}
n^{-\frac{1}{2}-\delta}\kappa_{\eta}^{-1}(\delta)S_{\left[  nr\right]
}\Rightarrow\frac{1}{\Gamma(1+\delta)}\int_{0}^{r}(r-s)^{\delta}%
d\mathbf{w}(s),\text{ \ \ when }-\frac{1}{2}<\delta<\frac{1}{2}, \tag{$3.6$}%
\end{equation}
and%
\begin{equation}
n^{-\frac{1}{2}-\delta}\left(  \log^{-1}n\right)  \kappa_{\eta}^{-1}(\tfrac
{1}{2})S_{\left[  nr\right]  }\Rightarrow\mathbf{w}_{0.5}(r),\text{ \ when
}\delta=0.5, \tag{$3.7$}%
\end{equation}
where $\mathbf{w}(r)$ is the standard Brownian motion on $\left[  0,1\right]
$ associated with the $\left(  u_{t}\right)  $ sequence and the symbols
$"\Rightarrow"$ and $"\overset{p}{\rightarrow}"$ denotes respectively weak
convergence and convergence in probability.

Since $d-d_{0}$ can always be decomposed as $d-d_{0}=m+\delta$, where $m\in%
\mathbb{N}
$ and $\delta\in]-0.5,0.5]$, the following result provides the asymptotic
distribution of the Dickey-Fuller normalized bias statistic $n\widehat{\rho
}_{n}=n\left(  \widehat{\phi}_{n}-1\right)  $ and the Dickey-Fuller
t-statistic, $t_{\widehat{\rho}_{n}}$, in the least square esimate of the
model $(3.3)$.

\begin{theorem}
\textit{Let }$\left\{  y_{t}\right\}  _{t=1}^{n}$\textit{\ be generated from
the DGP }$(2.1)$\textit{. If the regression model }$(3.3)$\textit{\ is fitted
to a sample of size }$n$\textit{\ then, as }$n\uparrow\infty$\textit{,}

\begin{enumerate}
\item $n\widehat{\rho}_{n}$\textit{\ satisfies}%
\begin{equation}
\widehat{\rho}_{n}=O_{p}(1)\text{ and }n\widehat{\rho}_{n}\overset
{p}{\rightarrow}-\infty,\text{ if }-1\leq d-d_{0}<-0.5, \tag{$3.8$}%
\end{equation}%
\begin{equation}
\widehat{\rho}_{n}=O_{p}(\log^{-1}n)\text{ and }n\widehat{\rho}_{n}\overset
{p}{\rightarrow}-\infty,\text{ if }d-d_{0}=-0.5, \tag{$3.9$}%
\end{equation}%
\begin{equation}
\widehat{\rho}_{n}=O_{p}(n^{-1-2\delta})\text{ and }n\widehat{\rho}%
_{n}\overset{p}{\rightarrow}-\infty,\text{ if }-0.5<d-d_{0}<0, \tag{$3.10$}%
\end{equation}%
\begin{equation}
\widehat{\rho}_{n}=O_{p}(n^{-1})\text{ and }n\widehat{\rho}_{n}\Rightarrow
\frac{\frac{1}{2}\left\{  \mathbf{w}^{2}(1)-1\right\}  }{\int_{0}%
^{1}\mathbf{w}^{2}(r)dr},\text{ if }d-d_{0}=0, \tag{$3.11$}%
\end{equation}%
\begin{equation}
\widehat{\rho}_{n}=O_{p}(n^{-1})\text{ and }n\widehat{\rho}_{n}\Rightarrow
\frac{\frac{1}{2}\mathbf{w}_{\delta,m+1}^{2}(1)}{\int_{0}^{1}\mathbf{w}%
_{\delta,m+1}^{2}(r)dr},\text{ if }d-d_{0}>0. \tag{$3.12$}%
\end{equation}

\item $t_{\widehat{\rho}_{n}}$\textit{\ is such that}%
\begin{equation}
t_{\widehat{\rho}_{n}}=O_{p}(n^{0.5})\text{ and }t_{\widehat{\rho}_{n}%
}\overset{p}{\rightarrow}-\infty,\text{ if }-1\leq d-d_{0}<-0.5, \tag{$3.13$}%
\end{equation}
\textit{ }%
\begin{equation}
t_{\widehat{\rho}_{n}}=O_{p}(n^{0.5}\log^{0.5}n)\text{ and }t_{\widehat{\rho
}_{n}}\overset{p}{\rightarrow}-\infty,\text{ if }d-d_{0}=-0.5, \tag{$3.14$}%
\end{equation}%
\begin{equation}
t_{\widehat{\rho}_{n}}=O_{p}(n^{-\delta})\text{ and }t_{\widehat{\rho}_{n}%
}\overset{p}{\rightarrow}-\infty,\text{ if }-\frac{1}{2}<d-d_{0}<0,
\tag{$3.15$}%
\end{equation}%
\begin{equation}
t_{\widehat{\rho}_{n}}=O_{p}(1)\text{ and }t_{\widehat{\rho}_{n}}%
\Rightarrow\frac{\frac{1}{2}\left\{  \mathbf{w}^{2}\left(  1\right)
-1\right\}  }{\left[  \int_{0}^{1}\mathbf{w}^{2}(r)dr\right]  ^{1/2}}\text{,
if }d-d_{0}=0, \tag{$3.16$}%
\end{equation}%
\begin{equation}
t_{\widehat{\rho}_{n}}=O_{p}(n^{\delta})\text{ and }t_{\widehat{\rho}_{n}%
}\overset{p}{\rightarrow}+\infty,\text{ if }0<d-d_{0}<0.5, \tag{$3.17$}%
\end{equation}%
\begin{equation}
t_{\widehat{\rho}_{n}}=O_{p}(n^{0.5})\text{ and }t_{\widehat{\rho}_{n}%
}\overset{p}{\rightarrow}+\infty,\text{ if }d-d_{0}\geq0.5. \tag{$3.18$}%
\end{equation}
\textit{where }$\mathbf{w}_{\delta,m}(r)$\textit{\ is the }$(m-1)-$%
\textit{fold integral of }$\mathbf{w}_{\delta}(r)$\textit{\ recursively
defined as }$\mathbf{w}_{\delta,m}(r)=\int_{0}^{r}\mathbf{w}_{\delta
,m-1}(s)ds$, \textit{with }$\mathbf{w}_{\delta,1}(r)=\mathbf{w}_{\delta}%
(r)$\textit{\ and }$\mathbf{w}(r)$\textit{\ denotes the standard Brownian
motion}.
\end{enumerate}
\end{theorem}

\begin{proof}
See Appendix.
\end{proof}

The later properties represent generalizations of those established by Sowell
$\left[  23\right]  $ for the cases $-\frac{1}{2}<d-1<0$, $d-1=0$ and
$0<d-1<\frac{1}{2}$. From $(3.8)$, $(3.9)$ and $\left(  3.10\right)  $, the
rate at which $\widehat{\rho}_{n}=\widehat{\phi}_{n}-1$ converges to zero
(i.e. $\widehat{\phi}_{n}$ converge to $1$ ) is slow for non-positive values
of $d-d_{0}$, and is particularly very slow for $-\frac{1}{2}<d-d_{0}%
<-\frac{1}{4}$. Moreover for $-\frac{1}{2}<d-d_{0}<0$, the limiting
distribution of $\widehat{\rho}_{n}$ has non-positive support and then
$\underset{n\rightarrow\infty}{\lim}P\left(  \widehat{\phi}_{n}<1\right)  =1$.
From $(3.11)$ and $(3.12)$, $\widehat{\rho}_{n}$ converges to zero at rate
$n$, when $d\geq d_{0}$. The rate of convergence $n$ is faster than the usual
standard rate $n^{\frac{1}{2}}$, when we deal with stationary $I(0)$
variables. Then, for $d-d_{0}\geq0$, the least squares estimate is
superconsistent. In other words, if a first order autoregression $(3.2)$ is
fitted to a series generated from an $ARFIMA(0,1+d-d_{0},0)$, where
$1+d-d_{0}$ is the order of integration of $\Delta^{-1+d_{0}}y_{t}$, then when
$d-d_{0}\geq0$, asymptotically, the $OLS$ estimator $\widehat{\phi}_{n}$ will
not exceed $1$ in probability. Figure $1$ and Figure $2$ below illustrate this
fact in an obvious way.

\begin{figure}[ptbh]
\caption{Relation between the order of integration, $d$, of the process
$y_{t}$ and the OLS estimator, $\widehat{\phi}_{n}$, in the regression model
$\Delta^{-1+d_{0}}y_{t}=\phi\Delta^{-1+d_{0}}y_{t-1}+\varepsilon_{t},$
($d_{0}$ fixed and $d$ varied)}
\begin{center}
\includegraphics[width=7cm]{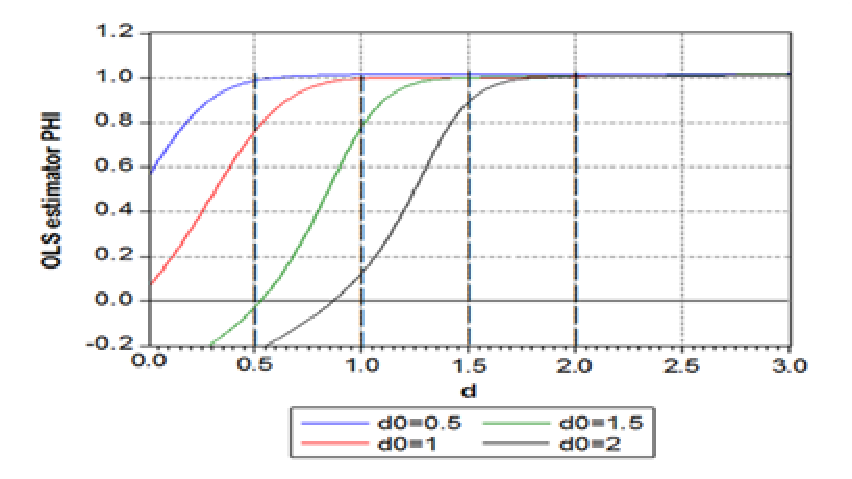}
\end{center}
\par
\end{figure}Figure $1$ shows that $\widehat{\phi}_{n}=1$ as long as
$d-d_{0}\geq0$, and $\widehat{\phi}_{n}<1$ as long as $d-d_{0}<0$, where
$\widehat{\phi}_{n}$ is the $OLS$ estimator in the autoregression model
$(3.2)$.

\begin{example}
For example, for $d_{0}=0.5$, we have$,$%
\[
\left\{
\begin{array}
[c]{c}%
d-0.5<0\text{ and }\widehat{\phi}_{n}<1\text{ for }0\leq d<0.5,\\
d-0.5\geq0\text{ and }\widehat{\phi}_{n}=1\text{ for }d\geq0.5,
\end{array}
\right.
\]
and for $d_{0}=2$, we have,%
\[
\left\{
\begin{array}
[c]{c}%
d-2<0\text{ and }\widehat{\phi}_{n}<1\text{ for }0\leq d<2,\\
d-2\geq0\text{ and }\widehat{\phi}_{n}=1\text{ for }d\geq2.
\end{array}
\right.
\]

\end{example}

\begin{figure}[ptbh]
\caption{Relation between the order of integration, $d$, of the process
$y_{t}$ and the $OLS$ estimator, $\widehat{\phi}_{n}$, in the regression model
$\Delta^{-1+d_{0}}y_{t}=\phi\Delta^{-1+d_{0}}y_{t-1}+\varepsilon_{t},$ ($d$
fixed and $d_{0}$ varied)}
\begin{center}
\includegraphics[width=7cm]{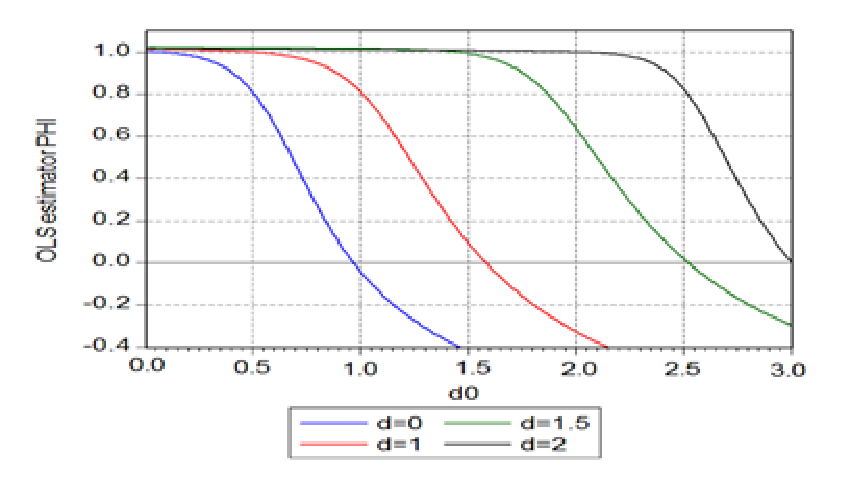}
\end{center}
\par
\end{figure}

Figure $2$ shows that as long as $d-d_{0}\geq0$, we have $\widehat{\phi}%
_{n}=1$, and $\widehat{\phi}<1$ whenever $d-d_{0}<0,$ where $\widehat{\phi
}_{n}$ is the $OLS$ estimator in the autoregression model $(3.2)$.

\begin{example}
For example, when $d=0.5$,%
\[
\left\{
\begin{array}
[c]{c}%
0.5-d_{0}<0\text{ and }\widehat{\phi}_{n}<1\text{ for }0\leq d_{0}<0.5,\\
0.5-d_{0}\geq0\text{ and }\widehat{\phi}_{n}=1\text{ for }d_{0}\geq0.5,
\end{array}
\right.
\]
and when $d=2,$%
\[
\left\{
\begin{array}
[c]{c}%
2-d_{0}<0\text{ and }\widehat{\phi}_{n}<1\text{ for }0\leq d_{0}<2,\\
2-d_{0}\geq0\text{ and }\widehat{\phi}_{n}=1\text{ for }d_{0}\geq2.
\end{array}
\right.
\]

\end{example}

Figure $1$ is made as follows: For a fixed sample $\left\{  u_{0}%
,\cdots,u_{1000}\right\}  $ generated from a Gaussian $i.i.d.(0,1)$ process,
samples of $ARFIMA(0,1+d-d_{0},0)$ processes were generated for $d$ varying
between $0$ and $3$ with step-size $0.01$ and $d_{0}$ fixed. Similarly, Figure
$2$ is made as follows. For a series $\left\{  u_{0},\cdots,u_{1000}\right\}
$ generated from a Gaussian $i.i.d.(0,1)$ process, samples from
$ARFIMA(0,1+d-d_{0},0)$ processes were generated for $d_{0}$ varying between
$0$ and $3$ with step-size $0.01$ for fixed $d$. For each series $\left\{
x_{t},t=1,\cdots,1000\right\}  $, a first order autoregression $(3.2)$ is
fitted and an estimate of $\phi$ is calculated. By plotting the estimate
$\widehat{\phi}_{n}$ against the fractional parameter $d$, one obtains Figure
$1$ and by plotting the parameter $\widehat{\phi}_{n}$ against the fractional
parameter $d_{0}$ one obtains Figure $2$. A general procedure for generating a
fractionally integrated series with length $n$ is to apply the formula
$x_{t}=\sum_{j=0}^{t-1}\frac{\Gamma\left(  d+1-d_{0}+j\right)  }{\Gamma\left(
d+1-d_{0}\right)  \Gamma\left(  j+1\right)  }u_{t-j}$ for $t=1,...,n$.

\begin{remark}
\item By fixing the parameter $d_{0}$ and varying the parameter $d$, we
increase the order of integration of $x_{t}$, and by varying the parameter
$d_{0}$ and fixing the parameter $d$ we decrease the order of integration of
$x_{t}$.
\end{remark}

The relationships between $\widehat{\phi}_{n}$ and $d$ and between
$\widehat{\phi}_{n}$ and $d_{0}$, highlighted by the results $(3.8)$-$(3.12)$
and illustrated by Figures $1$ and $2$, suggest that when we deal with testing
the degree of fractional integration, we have%
\[%
\begin{array}
[c]{ccc}%
H_{0}:d\geq d_{0}\Rightarrow H_{0}:\phi=1 & \text{and} & H_{1}:d<d_{0}%
\Rightarrow H_{1}:\phi<1.
\end{array}
\]

Like $DGM$, we call the test which is based on the hypotheses $(3.1)$ and the
auxiliary regression model $(3.2)$, or equivalently $(3.3)$ as "Fractional
Dickey-Fuller test ($F$-$DF$ test in short).

Another important property highlighted by Theorem $3.1$ is that the tests are
invariant to the original value of $d$, so the asymptotic properties only
depend on $d-d_{0}$. For example, we have used several series with sample size
$10000$ to estimate the densities (following Sowell $\left[  23\right]  $) of
$n\widehat{\rho}_{n}$ and $t_{\widehat{\rho}_{n}}$ under $d-d_{0}=0$. The
estimated densities are presented in Figures $3$ and $4$ above.

\begin{figure}[ptbh]
\caption{Kernel density estimate of $n\widehat{\rho}$ statistic under
$H_{0}:d=d_{0}$ using $1000$ samples of size $n=250$. $\left(  \text{\ }%
n\widehat{\rho}_{n}\Rightarrow\frac{\frac{1}{2}\left\{  \mathbf{w}%
^{2}(1)-1\right\}  }{\int_{0}^{1}\mathbf{w}^{2}(r)dr}\right)  $}
\begin{center}
\includegraphics[width=10cm]{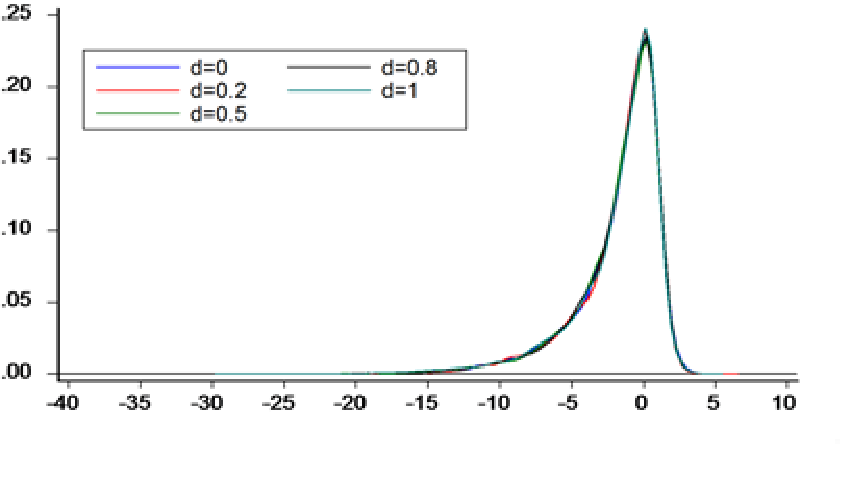}
\end{center}
\par
\end{figure}

\begin{figure}[ptbh]
\caption{Kernel density estimate of $t_{\widehat{\rho}\text{ }}$statistic
under $H_{0}:d=d_{0}$ using $1000$ samples of size $n=250$. \ $\left(
t_{\widehat{\rho}_{n}}\Rightarrow\frac{\frac{1}{2}\left\{  \mathbf{w}%
^{2}\left(  1\right)  -1\right\}  }{\left[  \int_{0}^{1}\mathbf{w}%
^{2}(r)dr\right]  ^{1/2}}\right)  $}
\begin{center}
\includegraphics[width=10cm]{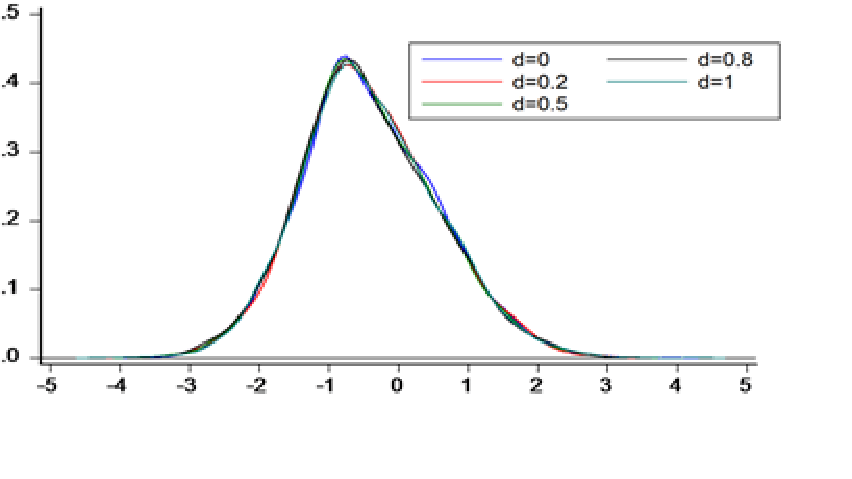}
\end{center}
\par
\end{figure}


In figures $3$ and $4$, for each one of the statistics $n\widehat{\rho}_{n}$
and $t_{\widehat{\rho}_{n}}$ under $d=d_{0}$ with $n=250$, the estimated
densities for different values of $d$ are represented on the same graph.
Figures $3$ and $4$ show that by fitting the regression model $(3.3)$ to the
sample generated from $(2.1)$, one obtains the same distribution as those used
by Dickey-Fuller $\left[  7\right]  $. In other words, as shown below, the
proposed test, which is based on the regression model $(3.2)$ (or equivalently
$(3.3)$) and the composite hypotheses $(3.1)$, can be understood and
implemented exactly as the simple Dickey-Fuller test for unit root by using
the usual statistical tables of the conventional statistics $n\widehat{\rho
}_{n}$ and ${\Large t}_{\widehat{{\Large \rho}}_{{\Large n}}}$.



\section{Power and size of the $F$-$DF$ test}

\subsection{Some theoretical aspects}

Let $Z_{1}={\footnotesize n(}\widehat{\phi}_{n}{\footnotesize -1)=n}%
\widehat{\rho}_{n}${\footnotesize \ \ \ }and\ \ \ $Z_{2}=t_{\widehat{\rho}%
_{n}}.$ For a composite hypothesis, the parameter space $\Omega=%
\mathbb{R}
$ is divided into disjoint regions, $\Omega_{0}=\left[  d_{0},+\infty\right[
$ and $\Omega_{1}=\left]  -\infty,d_{0}\right[  $. The test is written%
\begin{equation}
H_{0}:d\in\Omega_{0}\text{ \ \ against \ \ \ }H_{1}:d\in\Omega_{1}.
\tag{$4.1$}%
\end{equation}
For a series generated from $(2.1)$ with sample size $n$ we introduce two
nonrandomized test defined by a function $\mathbf{\Psi}_{i,n}$, $i=1,2$ on the
sample space\ of the observations $Z_{i}$, $i=1,2$ with critical regions
$C_{i}$, $i=1,2$. The test $\mathbf{\Psi}_{i,n}$ for a given rejection region
$C_{i}$ is%
\begin{equation}
\mathbf{\Psi}_{i,n}\left(  z_{i}\right)  =\left\{
\begin{array}
[c]{c}%
1\text{ if }z_{i}\in C_{i},\\
0\text{ if }z_{i}\notin C_{i}\text{.}%
\end{array}
\right.  \tag{$4.2$}%
\end{equation}
The power of the test is defined by the function%
\[
\Pi_{\mathbf{\Psi}_{i,n}}(d)=P\left[  Z_{i}\in C_{i}\mid d\right]  .
\]
$\Pi_{\mathbf{\Psi}_{i,n}}(d)$ measures the probability of rejecting the null
hypothesis for a given $d$ and rejection region $C_{i}$. The ideal test
function has%
\[%
\begin{array}
[c]{ccc}%
\Pi_{\mathbf{\Psi}_{i,n}}(d)\approx0\text{ for all }d\in\Omega_{0} & \text{and
} & \Pi_{\mathbf{\Psi}_{i,n}}(d)\approx1\text{ for all }d\in\Omega_{1},
\end{array}
\]
and the test function yields the correct decision with probability nearly $1$.
The type $I$ and type $II$ errors can be summarized in the power function
$\Pi_{\mathbf{\Psi}_{i,n}}(d)$. For $d\in\Omega_{0},$%

\[
\Pi_{\mathbf{\Psi}_{i,n}}(d)\text{ is the probability of making a type
}I\text{ error (size of the test),}
\]
and for $d\in\Omega_{1}$,%
\[
1-\Pi_{\mathbf{\Psi}_{i,n}}(d)\text{ is the probability of making a type
}II\text{ error.}
\]
For the alternative hypothesis $H_{1}:d<d_{0}$, we consider the one sided
critical regions of the form%
\begin{equation}
C_{i}=\left\{  Z_{i}<c_{n,i}\left(  \alpha\right)  \right\}  , \tag{$4.3$}%
\end{equation}
where $\alpha$ is the level of the test and $c_{n,i}\left(  \alpha\right)  $
the critical points. The level $\alpha$ of the test $\mathbf{\Psi}_{i,n}$ is
given by%
\[
{\Large \alpha}=\underset{{\Large d\in\Omega}_{0}}{{\Large Sup}}{\Large \Pi
}_{{\Large \Psi}_{{\Large i,n}}}{\Large (d).}
\]
It measures the maximum probability of rejecting the null hypothesis when it
is true. For the statistic $Z_{1}$, the figures $5$ and $6$, where
$c=c_{n,1}\left(  \alpha\right)  $, show clearly that the supremum occurs at
$d=d_{0}$,%
\[
P\left[  Z_{1}\in C_{1}\mid d\in\Omega_{0}\right]  {\Large \leq}P\left[
Z_{1}\in C_{1}\mid d=d_{0}\right]  .
\]
For the statistic $Z_{2}$, the figure $7$ shows also clearly that%
\[
P\left[  Z_{2}\in C_{2}\mid d\in\Omega_{0}\right]  {\Large \leq}P\left[
Z_{2}\in C_{2}\mid d=d_{0}\right]  ,
\]
because the asymptotic distribution of $Z_{2}$ is well defined, in all the
real line for $d=d_{0}$ and diverge to $+\infty$ for $d>d_{0}$. The figures
$5$, $6$ and $7$ show that the level of the test is%
\begin{equation}
{\Large \alpha}=\underset{{\Large d\in\Omega}_{0}}{{\Large Sup}}{\Large \Pi
}_{{\Large \Psi}_{{\Large i,n}}}{\Large (d)=\Pi}_{{\Large \Psi}_{{\Large i,n}%
}}{\Large (d}_{0}{\Large ).} \tag{$4.4$}%
\end{equation}
Another technical argument that supports $(4.4)$, is that the asymptotic
distributions of $Z_{1}$ and $n^{-\delta}Z_{2}$, $\delta>0$, have positive
support for $d>d_{0}$ and well defined in all the real line for $d=d_{0}$.
Consequently, the critical points $c_{n,i}\left(  \alpha\right)  ,i=1,2$ are
those used in the simple Dickey-Fuller test (i.e. without trend and
intercept). Note that all the critical values $c_{n,i}\left(  \alpha\right)  $
are less than zero for $\alpha=1\%$, $5\%$ and $1\%$. As indicated by $(4.4)$,
a test has level $\alpha$ if its size is less than $\alpha.$%
\begin{figure}[ptb]
\begin{center}
\includegraphics[width=13cm]{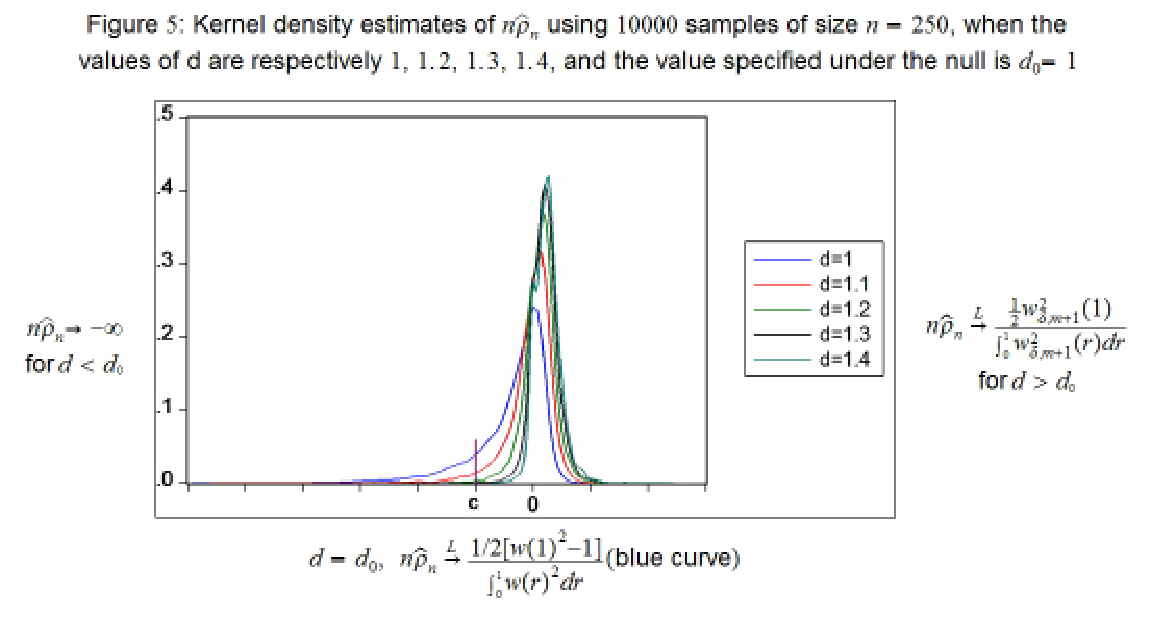}
\end{center}
\end{figure}

\begin{center}
\begin{figure}[ptb]
\begin{center}
\includegraphics[width=13cm]{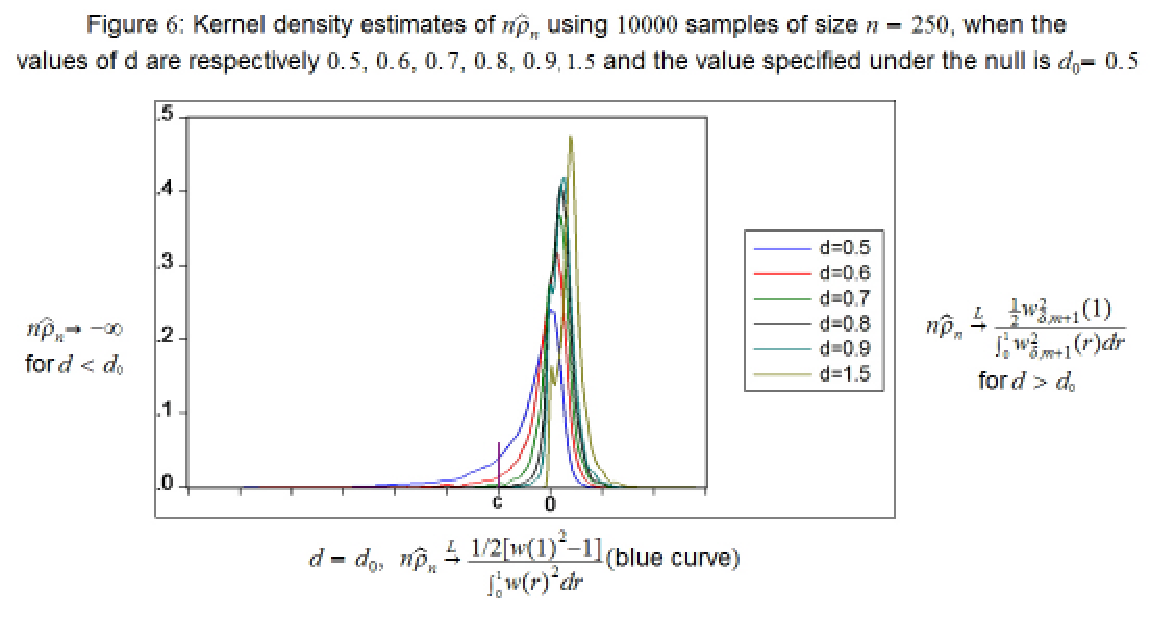}
\end{center}
\end{figure}

\begin{figure}[ptb]
\begin{center}
\includegraphics[width=13cm]{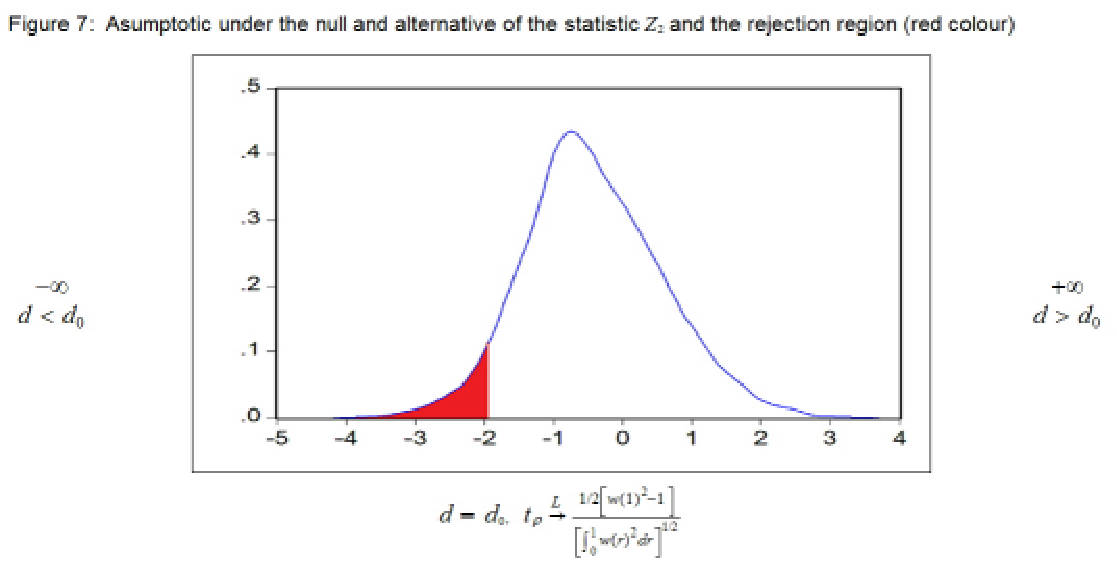}
\end{center}
\end{figure}
\end{center}

\begin{theorem}
For a given level $\alpha$, a sequence of test $\left\{  {\Large \Psi
}_{{\Large i,n}}\right\}  $, $i=1,2$, defined by (4.1), with critical region
(4.2) is consistent i.e.%
\begin{align*}
\underset{n\rightarrow\infty}{Lim}{\Large \Pi}_{{\Large \Psi}_{{\Large i,n}}%
}{\Large (d)}  &  {\Large =0,}\text{ for }d-d_{0}\geq0,\\
\underset{n\rightarrow\infty}{Lim}{\Large \Pi}_{{\Large \Psi}_{{\Large i,n}}%
}{\Large (d)}  &  {\Large =0,}\text{ for }-0.5<d-d_{0}<0.
\end{align*}

\end{theorem}

\begin{proof}
First, we consider the statistic $Z_{1}$. For $-0.5<d-d_{0}<0$ (i.e.
$-0.5<\delta<0$), since from result, $(3.10)$, the asymptotic distribution of
$n^{2\delta}Z_{1}$ has non-positive support and $\underset{n\rightarrow\infty
}{Lim}n^{2\delta}c_{n,1}\left(  \alpha\right)  =0$, then we have
\begin{align*}
\underset{n\rightarrow\infty}{Lim}P\left[  Z_{1}<c_{n,1}\left(  \alpha\right)
\mid d<d_{0}\right]   &  =\underset{n\rightarrow\infty}{Lim}P\left[
n^{2\delta}Z_{1}<n^{2\delta}c_{n,1}\left(  \alpha\right)  \right]  ,\\
&  =1.
\end{align*}
We have, also, $\underset{n\rightarrow\infty}{Lim}P\left[  Z_{1}%
<c_{n,1}\left(  \alpha\right)  \mid d>d_{0}\right]  =0$, because $Z_{1}$ has a
positive support for $d>d_{0}$, $c_{n,1}\left(  \alpha\right)  <0$ and
$\underset{n\rightarrow\infty}{Lim}n^{2\delta}c_{n,1}\left(  \alpha\right)
=0.$ Now consider the statistic $Z_{2}$. For For $-0.5<d-d_{0}<0$ (i.e.
$-0.5<\delta<0$), using the same arguments as above we have,%
\begin{align*}
\underset{n\rightarrow\infty}{Lim}P\left[  Z_{2}<c_{n,2}\left(  \alpha\right)
\mid d<d_{0}\right]   &  =\underset{n\rightarrow\infty}{Lim}P\left[
n^{\delta}Z_{2}<n^{\delta}c_{n,2}\left(  \alpha\right)  \right]  ,\\
&  =1.
\end{align*}
For $d>d_{0}$, with $\gamma=\delta$, when $0<\delta=d-d_{0}<0.5$ and
$\gamma=0.5$ when $d-d_{0}\geq0.5$, we have%
\[
\underset{n\rightarrow\infty}{Lim}P\left[  n^{-\gamma}Z_{2}<n^{-\gamma}%
c_{2,n}(\alpha)\mid d>d_{0}\right]  =0.
\]

\end{proof}

\subsection{Simulation study}

In this subsection, through a Monte Carlo study, we show that the proposed
$FD$-$F$ test performs very well in terms of power and size when we use the $t
$ statistic. To investigate the size and power of the $F$-$DF$ test, $10000$
samples of $FI(d)$ Gaussian processes $(2.1)$ are generated and the regression
model $(3.3)$ is used to estimate $t$. The sample-sizes considered are $n=50$
and $n=250$. Three values of $d$ are used$:0;$ $0.5$; $1$. For each value, we
specify the various values for $d_{0}$. Letting $S_{d}(d_{0})$ be the set of
values of $d_{0}$ for a given value of $d$, the sets which will be used for
the three values of $d$ are respectively%
\begin{align*}
S_{0}(d_{0})  &  =\left\{  -0.4\text{; }-0.3\text{; }-0.2\text{; }-0.1\text{;
}\mathbf{0}\text{; }0.1\text{; }0.2\text{; }0.3\text{; }0.4\right\}  ,\\
S_{0.5}(d_{0})  &  =\left\{  0\text{; }0.1\text{; }0.2\text{; }0.3\text{;
}0.4\text{; }\mathbf{0.5}\text{; }0.6\text{; }0.7\text{; }0.8\text{;
}0.9\right\}  ,\\
S_{1}(d_{0})  &  =\left\{  0.5;\text{ }0.6\text{; }0.7\text{; }0.8\text{;
}0.9\text{; }\mathbf{1}\text{; }1.1\text{; }1.2\text{; }1.3\text{;
}1.4\right\}  .
\end{align*}
Table $2$ and\ Table $3$ give simulation results on the size of the test,
(i.e. when $d-d_{0}=\delta\geq0$), where it may be easily seen that the
$F$-$DF$ test, based on the auxiliary regression model $(3.3)$, has good
performances in terms of size since%
\[
P\left(  t<c_{n}(\alpha\right)  \mid\delta>0)\leq P\left(  t<c_{n}%
(\alpha\right)  \mid\delta=0)\text{ and }P\left(  t<c_{n}(\alpha\right)
\mid\delta=0)\approx\alpha.
\]
Table $2$ and Table $3$ provide, also, the simulation results on the power of
the test (i.e. when $d-d_{0}=\delta<0$)%
\[
P\left(  t<c_{n}(\alpha\right)  \mid\delta<0)\geq P\left(  t<c_{n}%
(\alpha\right)  \mid\delta=0).
\]
In this case, there are some conclusions to be drawn from it. First, the power
of the $F$-$DF$ test increases with the increase of sample size and
$\delta=d-d_{0}$. For example, for $\alpha=5\%$, $d=1$ and $\delta=-0.1$, the
power is $12.36\%$ for $n=50$ and $20.76\%$ for $n=250$. When $\alpha=5\%$,
$d=1$ and $\delta=-0.3$, the power is $48.5\%$ for $n=50$ and $86.05\%$ for
$n=250$. Second, as shown in table $3$, for $n=250,$ the power of the $F$-$DF$
test is below $50\%$ for ($\delta=-0.1$) and for ($\alpha=1\%$, $\delta=-0.2
$). Third, for given $n$, $\alpha$ and $\delta$, the power for $d=0$, $d=0.5$
and $d=1$ are approximately similar because the asymptotic under the
alternative does not depend on $d$ but only on $\delta=d-d_{0}$. Finally,
another important property showed by the table $2$ and $3$ is that the power
function satisfies%
\begin{align*}
{\Large \Pi}_{{\Large \Psi}_{{\Large i,n}}}{\Large (d)}  &  {\Large \leq
\alpha}\text{ \ \ \ if }d-d_{0}\geq0,\\
{\Large \Pi}_{{\Large \Psi}_{{\Large i,n}}}{\Large (d)}  &  {\Large \geq
\alpha}\text{ \ \ \ if }d-d_{0}<0.
\end{align*}
A test for which the power function satisfies the conditions above is said to
be unbiased.

\begin{figure}[ptbh]
\begin{center}
\includegraphics[width=10cm]{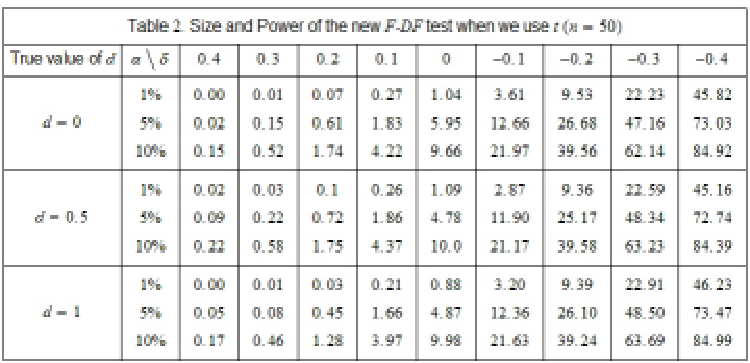}
\end{center}
\end{figure}

\begin{figure}[ptbh]
\begin{center}
\includegraphics[width=10cm]{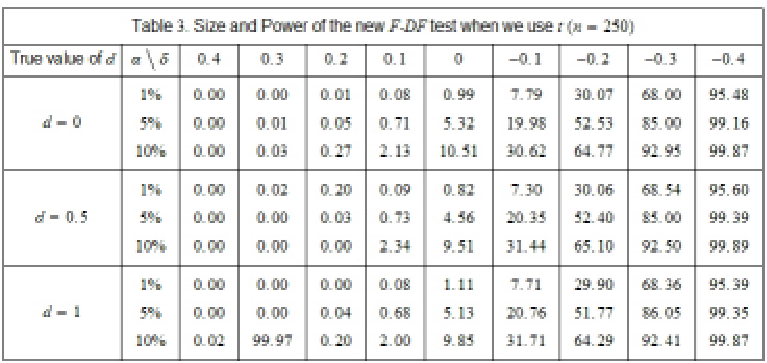}
\end{center}
\end{figure}

Similar results are obtained for the $n\widehat{\rho}_{n}$ statistic. Since
$n\widehat{\rho}_{n}$ has the non degenerate limit distribution, we choose to
give the simulation results in the form of estimation density (by kernel
methods, like Sowell $\left[  23\right]  $. The figures $5$ and $6$ above
summarize these results and support clearly those of Theorem $4.1$.\newpage

\section{Application to the Nelson-Plosser data}

For the sake of illustration, this section applies our $F$-$DF$ test to the
well-known Nelson-Plosser data. The starting date is $1860$ for the consumer
price index and industrial production, $1869$ for velocity, $1871$ for stock
prices, $1889$ for GNP deflator and money stock, $1890$ for employment and
unemployment rate, $1900$ for bond yield, real wages and wages, and $1909$ for
the nominal and real GNP and GNP per capita. The variables are expressed in
natural logarithms. All variables exhibit an upward trend with the exception
of velocity, which shows a strong downward trend and the unemployment rate
which tends to fluctuate around a constant level. The seminal empirical work
by Nelson and Plosser $\left[  18\right]  $ suggests that there is a strong
evidence for the unit root hypothesis for most macroeconomic time series data.
Two possible specifications for the data generating processes $(DGP)$ are then%
\begin{equation}
y_{t}=(1-L)^{-d}u_{t}, \tag{5.1}%
\end{equation}
and%
\begin{equation}
y_{t}=\alpha+(1-L)^{-d}u_{t}. \tag{5.2}%
\end{equation}
The theoretical framework provided in this paper does not allow us to use the
$DGP$ $(5.2)$ (see Appendix $2$). At this level, we only use the $DGP$
$(5.1)$. For the $DGP$ $(5.1)$, we test the null for several values of $d_{0}%
$, namely: $0$; $0.5$; $1$; $1.5$ and $2$ by using respectively the following
regression models,%
\begin{equation}
y_{t}=\rho\Delta^{-1}y_{t-1}+\varepsilon_{1,t}, \tag{Model (I)}%
\end{equation}%
\begin{equation}
\Delta^{0.5}y_{t}=\rho\Delta^{-0.5}y_{t-1}+\varepsilon_{2,t}, \tag{Model (II)}%
\end{equation}%
\begin{equation}
\Delta y_{t}=\rho y_{t-1}+\varepsilon_{3,t}, \tag{Model (III)}%
\end{equation}%
\begin{equation}
\Delta^{1.5}y_{t}=\rho\Delta^{0.5}y_{t-1}+\varepsilon_{4,t}, \tag{Model (IV)}%
\end{equation}%
\begin{equation}
\Delta^{2}y_{t}=\rho\Delta y_{t-1}+\varepsilon_{5,t}. \tag{Model (V)}%
\end{equation}
Note that, the sample sizes for all the $14$ U.S. macroeconomic Nelson-Plosser
series, used here, are between $n=80$ and $n=129$. Consequently, the decision
rules adopted for the testing problem $(2.3)$ are%
\begin{align*}
\text{reject }H_{0}\text{ if }Z_{1}  &  <-7.9,\\
\text{reject }H_{0}\text{ if }Z_{2}  &  <-1.95,
\end{align*}
where $Z_{1}$ and $Z_{2}$ are respectively the usual statistic $n\widehat
{\rho}$ and $\dfrac{\widehat{\rho}}{\sigma_{\widehat{\rho}}}$, and where
$\left(  -7.9,-1.95\right)  $ are the corresponding critical values at level
$\alpha=5\%$, obtained from the usual statistical tables of Dickey-Fuller
$\left[  5\right]  $. The results of the decision rules shown in Table $5$
suggest that:

\begin{itemize}
\item for model $(I)$, all series are found to be integrated with order
$d\geq0$,

\item for model $(II)$, all series are found to be integrated with order
$d\geq0.5$,

\item for model $(III)$, all series are found to be integrated with order
$d\geq1$,

\item for model $(IV)$, all series are found to be integrated with order
$d<1.5$, except the \textbf{Industrial production} and \textbf{Money stock }series.

\item for model $(V)$, all series are found to be integrated with order $d<2$.
\end{itemize}

In summary, it may be concluded from Table $5$ that, following our test, all
the macroeconomic variables are $d$-integrated with $1\leq d<1.5$, except for
the Industrial production and Money stock whose order of integration is
between $1.5$ and $2$, i.e. $1.5\leq d<2$. \begin{figure}[tbh]
\begin{center}
\includegraphics[width=8cm]{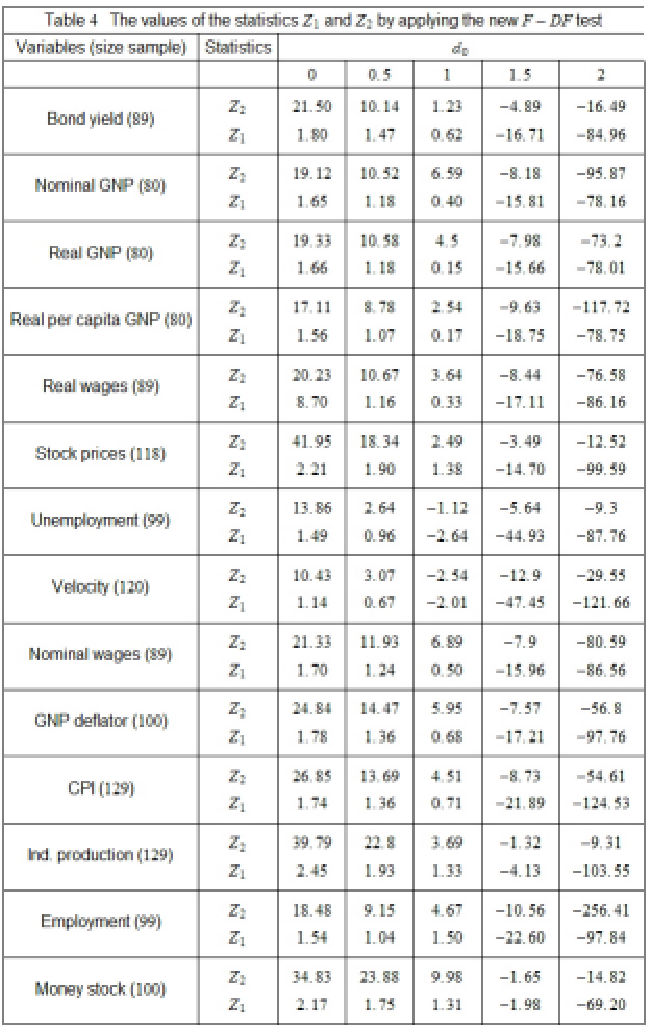}
\end{center}
\end{figure}

Note that the $F$-$DF$ test was done assuming that the empirical variables are
derived from data generating process $ARFIMA(0,d,0)$. A more general study is
needed to achieve adequate conclusions about the integration order for the
Nelson-Plosser Data, by considering more general data generating process
$ARFIMA(p,d,q)$ and also by incorporating non zero drift and time trend in
data generating process $(2.1)$ while using a suitable auxiliary regression
model. \newpage

\section{Concluding remarks and discussion.}

In the Dickey-Fuller paper, the parameter $d$, without restrict the
generality, can have only two values $d=1$ or $d=0$. To test%
\begin{equation}
H_{0}:d=1\text{ \ \ against \ \ }H_{1}:d=0, \tag{Standard test}%
\end{equation}
in the simple case, Dickey and Fuller use the regression model%
\begin{equation}
y_{t}=\phi y_{t-1}+\varepsilon_{t}. \tag{Standard regression model}%
\end{equation}

Since Anderson $\left[  2\right]  $, White $(\left[  25\right]  ,\left[
26\right]  )$ developed the statistical theory on the first order
autoregressive process with the autoregressive parameter equal to $1$(i.e.
$d=1$) and greater than one (explosive process), Box and Jenkins $\left[
5\right]  $ formalized the analysis of time series, and Nelson and Plosser
$\left[  19\right]  $ argue that the most macroeconomic series have unit
roots. The unit root test has been an important topic on the econometric literatures.

Phillips $\left[  20\right]  $ show that under the null hypothesis (i.e.
$H_{0}:d=1$) that the asymptotic distributions of $n\left(  \widehat{\phi
}-1\right)  $ and $t_{\widehat{\rho}}$ are respectively,

\begin{center}%
\begin{tabular}
[c]{|c|}\hline%
\begin{tabular}
[c]{ccc}%
$\frac{0.5\left\{  W^{2}(1)-1\right\}  }{\int_{0}^{1}W^{2}(r)dr}$ & and &
$\frac{0.5\left\{  W^{2}(1)-1\right\}  }{\left[  \int_{0}^{1}W^{2}%
(r)dr\right]  ^{0.5}}$\\
\multicolumn{3}{c}{Usual asymptotic distributions for DF test}%
\end{tabular}
\\\hline
\end{tabular}

\end{center}

where $W(\cdot)$ is the standard Brownien motion. These later asymptotic
distributions has been tabulated, the tabulated values are used to perform the
standard test.

In fractional case the parameter $d$ can have an infinite values, for example
$d$ can have an infinite number of values $\cdots-0.5$; $-0.1$; $-0.3 $;
$-0.4$; $-0.5$; $0;$ $0.1$; $0.2$; $0.3$; $0.4$; $0.5$;$\cdots$. For the
fractional case, the standard regression model can be used only for testing
the hypothesis $H_{0}:d=1$. In our paper, the question is

\begin{center}%
\begin{tabular}
[c]{|c|}\hline
$\text{\textbf{How to extend the standard framework above to take into
account}}$\\
$\text{\textbf{the fractional case?}}$\\\hline
\end{tabular}

\end{center}

Such extension has already been discussed by Dolado et Al $[10]$. Dolado et al
$\left[  10\right]  $ propose to test%
\[
H_{0}:d=d_{0}\text{ \ \ against \ \ }H_{1}:d=d_{1},
\]
by using the auxiliary regression model%
\begin{equation}
\Delta^{d_{0}}y_{t}=\rho\Delta^{d_{1}}y_{t-1}+\varepsilon_{t}. \tag{DGM
regression model}%
\end{equation}
In this paper, we show in the first step, that the $DGM$ approach is not the
best and adequate way to extend the Dickey-Fuller test by taking into account
the fractional case, because the $DGM$ regression model is based on the null
and the alternative (i.e. $d_{0}$ and $d_{1}$).

In the second step, we provide how to extend adequately the standard
Dickey-Fuller test $\left[  7\right]  $ by taking into account the fractional
case. In fact, in our approach, the question is

\begin{center}%
\begin{tabular}
[c]{|c|}\hline
$\text{\textbf{How to extend the standard framework below to take into
account}}$\\
$\text{\textbf{the fractional case by using the usual asymptotic
distribution}}$\\
$\frac{0.5\left\{  W^{2}(1)-1\right\}  }{\int_{0}^{1}W^{2}(r)dr}\text{ and
}\frac{0.5\left\{  W^{2}(1)-1\right\}  }{\left[  \int_{0}^{1}W^{2}%
(r)dr\right]  ^{0.5}}\text{\textbf{?}}$\\\hline
\end{tabular}

\end{center}

A correct answer to this question can be very useful in practice. The answer
we give to this question is based on four points:

\begin{enumerate}
\item Using the composite hypothesis $H_{0}:d\geq d_{0}.$

\item If $y_{t}\leadsto I(d_{0})$ than $(1-L)^{-1+d_{0}}y_{t}\leadsto I(1).$

\item Testing the composite null hypothesis is based upon testing the
statistical significance of the coefficient $\phi$ (or $\rho=\phi-1$) in the
regression model $\Delta^{-1+d_{0}}y_{t}=\phi\Delta^{-1+d_{0}}y_{t-1}%
+\varepsilon_{t}.$

\item The level of the test $\alpha=Sup_{d\geq d_{0}}P(reject$ $H_{0}%
)=P(reject$ $H_{0}|d=d_{0}).$
\end{enumerate}

Our test is based on a composite null hypothesis, $H_{0}:d\geq d_{0}$, this
choice was not done arbitrarily. This choice was made based on the results of
the asymptotic theory given in the theorem $4$. To use our test, we recommend
to follow the following steps:

\begin{enumerate}
\item Estimate the parameter $\rho$ in the regression model $\Delta^{d_{0}%
}y_{t}=\rho\Delta^{-1+d_{0}}y_{t-1}+\varepsilon_{t}.$ This regression provides
a more flexible and unified framework to test the null for different values of
$d_{0}$ while using the same critical value.

\item The null hypothesis is rejected if $Z_{i}<c_{i}(\alpha)$, where
$Z_{1}=t_{\widehat{\rho}_{n}}$ and $Z_{2}=n\widehat{\rho}_{n}$. The level of
the test can be approximated by its asymptotic value: $\alpha=Sup_{d\geq
d_{0}}P(Z_{i}<c_{i}(\alpha))=P[Z_{i}<c_{i}(\alpha)/d=d_{0}]\ $($i=1,2$).

\item The critical values $c_{i}(\alpha)$ ($i=1,2$)\ can be chosen so as to
achieve a predetermined size by using the usual Dickey-Fuller statistical tables.
\end{enumerate}

Finally, some remarks are in order:

- To implement our test we do not need to estimate the parameter $d$.

- We have referred to our test as the Fractional Dickey-Fuller ($F$-$DF$)
test. A similar designation, $F$-$DF$, has been adopted by Dolado et
\textit{al} $\left[  10\right]  $ for their test.

- Regarding the Dickey-Pantula test, both the upward and downward procedures
are still valid in our fractional case (see Dickey and Pantula $\left[
8\right]  $). Moreover, by sequentially repeating the test in upward or in the
downward senses, we can cover the value of $d$ at the desired accuracy.

- The empirical study on the Nelson-Plosser Data is only made to illustrate
the $F$-$DF$ test.

- In this article we have not discussed the situation when there is an
additional short memory component in the series, like the $AR$ or $MA$. Also,
the situation when there is a non-zero drift or a time trend in data
generating process may be investigated. In fact, the proposed $F$-$DF$ test
may be easily generalized to such situations. Here, we give just an indication
when $y_{t}\sim ARFIMA(p,d,0)$
\[
A(L)\Delta^{d}y_{t}=u_{t}\text{,}
\]
where $A(L)=\sum_{j=0}^{p}\alpha_{j}L^{j}$, $L$ is the backward shift
operator, $\alpha_{0}=1$, the roots of $A(z)=0$ are outside the unit circle
and $u_{t}$ is defined as above. Then the fractional augmented Dickey-Fuller
test, for the null hypothesis $d\geq d_{0}$, would be based on the regression
model%
\[
\Delta^{d_{0}}y_{t}=\rho\Delta^{-1+d_{0}}y_{t-1}+\sum_{j=0}^{p}\alpha
_{j}\Delta^{d_{0}}y_{t-j}+\varepsilon_{t}.
\]
Further research is currently being undertaken toward generalizing the
$F$-$DF$ testing approach, along similar directions as the $DF$ test has been
extended in the unit root literature accounting for time series which may
exhibit a trending behavior and for general $ARFIMA$ case.

\textbf{Appendix 1: Proof of Theorem 1}

By denoting $\Delta^{-1+d_{0}}y_{t}=x_{t}$, the $OLS$ estimator of $\rho$ and
its $t$-ratio for the auxiliary regression model $(3.3)$, are given by the
usual squares expressions%
\[
\widehat{\rho}_{n}=\frac{\sum_{t=1}^{n}\left(  \Delta x_{t}\right)  \left(
x_{t-1}\right)  }{\sum_{t=1}^{n}\left(  x_{t-1}\right)  ^{2}},\text{
\ \ \ \ \ \ \ \ }t_{\widehat{\rho}_{n}}=\frac{\sum_{t=1}^{n}\left(  \Delta
x_{t}\right)  \left(  x_{t-1}\right)  }{\left\{  \widehat{\sigma}_{n}^{2}%
\sum_{t=1}^{n}\left(  x_{t-1}\right)  ^{2}\right\}  ^{1/2}}%
\]
where the variance of the residuals, $\widehat{\sigma}_{n}^{2}$ is given by
$\widehat{\sigma}_{n}^{2}=n^{-1}\sum_{t=1}^{n}\left(  \Delta x_{t}%
-\widehat{\rho}_{n}x_{t-1}\right)  ^{2}$. Note that, $x_{t-1}\sim
FI(1+d-d_{0})$ and $\Delta x_{t}\sim FI(d-d_{0})$. Since $x_{t}$ is stationary
fractionally integrated process for $d-d_{0}\in\left[  -1;-0.5\right[  $ and
nonstationary fractional integrated process for $d-d_{0}\in\left[
-0.5;+\infty\right[  $, we divide our proof into two parts

\begin{enumerate}
\item \underline{\textbf{Part} $\mathbf{1:d-d}_{0}\mathbf{\in}\left[
-1;-0.5\right[  $}
\end{enumerate}

When $-1\leq d-d_{0}<-0.5$, given that $x_{t}$ is stationary fractionally
integrated, of order $\delta\in\left[  0;0.5\right[  $ and ergodic process,
then%
\[
\widehat{\rho}_{n}=\frac{\sum_{t=1}^{n}\left(  \Delta x_{t}\right)  \left(
x_{t-1}\right)  }{\sum_{t=1}^{n}\left(  x_{t-1}\right)  ^{2}}=\frac{n^{-1}%
\sum_{t=1}^{n}x_{t}x_{t-1}}{n^{-1}\sum_{t=1}^{n}x_{t-1}^{2}}-1\overset
{P}{\rightarrow}\frac{E\left(  x_{t}x_{t-1}\right)  }{E\left(  x_{t-1}%
^{2}\right)  }-1.
\]
Therefore, given that%
\begin{equation}
E\left(  x_{t}x_{t-j}\right)  =\sigma_{u}^{2}\frac{\Gamma(j+\delta
)\Gamma(1-2\delta)}{\Gamma(j+1-\delta)\Gamma(1-\delta)\Gamma(\delta)},\text{
\ }j\geq0, \tag{$A0$}%
\end{equation}
(see, Hosking $[12]$) and the recursive identity $\Gamma(1+z)=z\Gamma(z)$, it
follows that%
\[
\widehat{\rho}_{n}\overset{P}{\rightarrow}\frac{E\left(  x_{t}x_{t-1}\right)
}{E\left(  x_{t-1}^{2}\right)  }=-\frac{1}{1-\delta},
\]
which, in turn, given that $\delta\in\left[  0;0.5\right[  $, entails that
$\widehat{\rho}_{n}\in\left[  -2;-1\right[  $. Consequently, $n\widehat{\rho
}_{n}\overset{P}{\rightarrow}-\infty$. With respect to the $t$-test (i.e.
$t_{\widehat{\rho}_{n}}$), it is straightforward to prove that%
\[
n^{-1}\widehat{\sigma}_{n}^{2}\sum_{t=1}^{n}\left(  x_{t-1}\right)
^{2}=n^{-2}\left(  \left(  \sum_{t=1}^{n}x_{t}^{2}\right)  \left(  \sum
_{t=1}^{n}x_{t-1}^{2}\right)  -\left(  \sum_{t=1}^{n}x_{t}x_{t-1}\right)
^{2}\right)  ,
\]
and then, by using $(A_{0})$ and the recursive identity $\Gamma(1+z)=z\Gamma
(z)$, it follows that%
\[
\left(  n^{-1}\widehat{\sigma}_{n}^{2}\sum_{t=1}^{n}\left(  x_{t-1}\right)
^{2}\right)  ^{\frac{1}{2}}\overset{p}{\rightarrow}\sigma_{u}^{2}\frac
{\Gamma(1-2\delta)}{\Gamma\left(  1-\delta\right)  ^{2}}\frac{\sqrt{1-2\delta
}}{1-\delta},
\]
and%
\[
n^{-1}\left(  \sum_{t=1}^{n}x_{t}x_{t-1}-\sum_{t=1}^{n}x_{t-1}^{2}\right)
\overset{p}{\rightarrow}\sigma_{u}^{2}\frac{\Gamma(1-2\delta)}{\Gamma\left(
1-\delta\right)  ^{2}}\frac{2\delta-1}{1-\delta}.
\]
Consequently, $n^{-1/2}t_{\widehat{\rho}_{n}}\rightarrow-\sqrt{1-2\delta}$ and
then $t_{\widehat{\rho}_{n}}\overset{p}{\rightarrow}-\infty.$

\begin{enumerate}
\item[2] \underline{\textbf{Part }$\mathbf{2:d-d}_{\mathbf{0}}\mathbf{\in
}\left[  -0.5;+\infty\right[  $}.
\end{enumerate}

For the $\sum_{t=1}^{n}\left(  x_{t-1}\right)  ^{2}$ term, it follows from
$(3.4)$, $(3.5)$, $(3.6)$, $(3.7)$ and the continuous mapping theorem. When
$d-d_{0}=-0.5,$%
\begin{equation}
\frac{1}{n\left(  \log n\right)  \kappa_{\eta}^{2}(\frac{1}{2})}\overset
{n}{\underset{t=1}{\sum}}\left(  x_{t-1}\right)  ^{2}\Rightarrow\int_{0}%
^{1}\mathbf{w}_{0.5}^{2}(r)dr. \tag{A1}%
\end{equation}
When $-\frac{1}{2}<d-d_{0}<0,$%
\begin{equation}
\frac{1}{n^{2+2\delta}\kappa_{\eta}^{2}(\delta)}\overset{n}{\underset
{t=1}{\sum}}\left(  x_{t-1}\right)  ^{2}\Rightarrow\int_{0}^{1}\mathbf{w}%
_{\delta}^{2}(r)dr\text{.} \tag{A2}%
\end{equation}
When $d-d_{0}=0,$%
\begin{equation}
\frac{1}{n^{2}\kappa_{\eta}^{2}(0)}\overset{n}{\underset{t=1}{\sum}}\left(
x_{t-1}\right)  ^{2}\Rightarrow\int_{0}^{1}\mathbf{w}^{2}(r)dr\text{.}
\tag{A3}%
\end{equation}
When $0<d-d_{0}<\frac{1}{2},$%
\begin{equation}
\frac{1}{n^{2+2\delta}\kappa_{\eta}^{2}(\delta)}\overset{n}{\underset
{t=1}{\sum}}\left(  x_{t-1}\right)  ^{2}\Rightarrow\int_{0}^{1}\mathbf{w}%
_{\delta}^{2}(r)dr\text{.} \tag{A4}%
\end{equation}
When $d-d_{0}=m+\frac{1}{2},$ $m\geq1,$%
\begin{equation}
\frac{1}{n^{2(m+\frac{3}{2})}\left(  \log n\right)  \kappa_{\eta}^{2}(\frac
{1}{2})}\overset{n}{\underset{t=1}{\sum}}\left(  x_{t-1}\right)
^{2}\Rightarrow\int_{0}^{1}\mathbf{w}_{0.5,m+1}^{2}(r)dr\text{.} \tag{A5}%
\end{equation}
When $d-d_{0}=m+\delta$, $m\geq1,$%
\begin{equation}
\frac{1}{n^{2(m+1+\delta)}\kappa_{\eta}^{2}(\delta)}\overset{n}{\underset
{t=1}{\sum}}\left(  x_{t-1}\right)  ^{2}\Rightarrow\int_{0}^{1}\mathbf{w}%
_{\delta,m+1}^{2}(r)dr\text{.} \tag{A6}%
\end{equation}
For the $\sum_{t=1}^{n}\left[  \Delta x_{t}\right]  \left[  x_{t-1}\right]  $
term, we have%
\[
\sum_{t=1}^{n}\left[  \Delta x_{t}\right]  \left[  x_{t-1}\right]  =\frac
{1}{2}\left(  \Delta^{-1+d_{0}}y_{n}\right)  ^{2}-\frac{1}{2}\overset
{n}{\underset{t=1}{\sum}}\left(  \Delta x_{t}\right)  ^{2}
\]
For the first term, it follows from (3.4), (3.5), (3.6), (3.7) and the
continuous mapping theorem. When $d-d_{0}=-0.5,$%
\begin{equation}
\frac{1}{2\left(  \log n\right)  \kappa_{\eta}^{2}(\frac{1}{2})}\left(
\Delta^{-1+d_{0}}y_{n}\right)  ^{2}\Rightarrow\frac{1}{2}\mathbf{w}_{0.5}%
^{2}(1). \tag{A7}%
\end{equation}
When $-0.5<d-d_{0}<0,$%
\begin{equation}
\frac{1}{2n^{1+2\delta}\kappa_{\eta}^{2}(\delta)}\left(  \Delta^{-1+d_{0}%
}y_{n}\right)  ^{2}\Rightarrow\frac{1}{2}\mathbf{w}_{\delta}^{2}(1). \tag{A8}%
\end{equation}
When $d-d_{0}=0,$%
\begin{equation}
\frac{1}{2n^{1}\kappa_{\eta}^{2}(0)}\left(  \Delta^{-1+d_{0}}y_{n}\right)
^{2}\Rightarrow\frac{1}{2}\mathbf{w}^{2}(1). \tag{A9}%
\end{equation}
When $0<d-d_{0}<0.5,$%
\begin{equation}
\frac{1}{2n^{1+2\delta}\kappa_{\eta}^{2}(\delta)}\left(  \Delta^{-1+d_{0}%
}y_{n}\right)  ^{2}\Rightarrow\frac{1}{2}\mathbf{w}_{\delta}^{2}(1). \tag{A10}%
\end{equation}
When $d-d_{0}=m+0.5,$ $m\geq1,$%
\begin{equation}
\frac{1}{2n^{2\left(  m+1\right)  }\left(  \log n\right)  \kappa_{\eta}%
^{2}(\frac{1}{2})}\left(  \Delta^{-1+d_{0}}y_{n}\right)  ^{2}\Rightarrow
\frac{1}{2}\mathbf{w}_{0.5,m+1}^{2}(1). \tag{A11}%
\end{equation}
when $d-d_{0}=m+\delta,$ $m\geq1,$%
\begin{equation}
\frac{1}{2n^{1+2\left(  m+\delta\right)  }\kappa_{\eta}^{2}(\frac{1}{2}%
)}\left(  \Delta^{-1+d_{0}}y_{n}\right)  ^{2}\Rightarrow\frac{1}{2}%
\mathbf{w}_{\delta,m+1}^{2}(1). \tag{A12}%
\end{equation}
For the second term, we have:\newline When $d-d_{0}=-0.5$, by using Lemma
$2.1$ of Ming Liu $(1998)$ result $2$%
\begin{equation}
-\frac{1}{2}\kappa_{\eta}^{-2}(\frac{1}{2})n^{-1}\overset{n}{\underset
{t=1}{\sum}}\left(  \Delta x_{t}\right)  ^{2}\overset{p}{\rightarrow}-\frac
{1}{2}\kappa_{\eta}^{-2}(\frac{1}{2})var(\Delta x_{t})=-1. \tag{A13}%
\end{equation}
When $-0.5<d-d_{0}<0$, by using $(3.4)$ and the ergodic theorem (note that
here $d-d_{0}=\delta$)%
\begin{equation}
-\frac{1}{2}\kappa_{\eta}^{-2}(\delta)n^{-1}\overset{n}{\underset{t=1}{\sum}%
}\left(  \Delta x_{t}\right)  ^{2}\overset{p}{\rightarrow}-\frac{1}{2}%
\kappa_{\eta}^{-2}(\delta)var(\Delta x_{t})=-(\frac{1}{2}+\delta)\frac
{\Gamma(1+\delta)}{\Gamma\left(  1-\delta\right)  }. \tag{A14}%
\end{equation}
When $d-d_{0}=0$, by using $(3.4)$ and the ergodic theorem%
\begin{equation}
-\frac{1}{2}\kappa_{\eta}^{-2}(0)n^{-1}\overset{n}{\underset{t=1}{\sum}%
}\left(  \Delta x_{t}\right)  ^{2}\overset{p}{\rightarrow}-\frac{1}{2}%
\kappa_{\eta}^{-2}(0)var(\Delta x_{t})=-\frac{1}{2}. \tag{A15}%
\end{equation}
When $0<d-d_{0}<0.5$, by using $(3.4)$ and the ergodic theorem (note that here
$d-d_{0}=\delta$)%
\begin{equation}
-\frac{1}{2}n^{-1}\overset{n}{\underset{t=1}{\sum}}\left(  \Delta
x_{t}\right)  ^{2}\overset{p}{\rightarrow}-\frac{1}{2}\kappa_{\eta}%
^{-2}(\delta)var(\Delta x_{t})=-(\frac{1}{2}+\delta)\frac{\Gamma(1+\delta
)}{\Gamma\left(  1-\delta\right)  }. \tag{A16}%
\end{equation}
When $d-d_{0}=m+0.5$, $m\geq1$, by using $(3.5)$, $(3.7)$ and the continuous
mapping theorem%
\begin{equation}
\frac{1}{n^{2(m+0.5)}\left(  \log n\right)  \kappa_{\eta}^{2}(0.5)}\overset
{n}{\underset{t=1}{\sum}}\left(  \Delta x_{t}\right)  ^{2}\Rightarrow\int
_{0}^{1}\mathbf{w}_{0.5,m}^{2}(r)dr. \tag{A17}%
\end{equation}
When $d-d_{0}=m+\delta$, $m\geq1$, by using $(3.4)$, $(3.6)$ and the
continuous mapping theorem%
\begin{equation}
\frac{1}{n^{2(m+\delta)}\kappa_{\eta}^{2}(\delta)}\overset{n}{\underset
{t=1}{\sum}}\left(  \Delta x_{t}\right)  ^{2}\Rightarrow\int_{0}^{1}%
\mathbf{w}_{\delta,m}^{2}(r)dr. \tag{A18}%
\end{equation}
\noindent Therefore, when $d-d_{0}=-0.5$, we have, using $(A7)$ and $(A13)$,%
\begin{equation}
n^{-1}\kappa_{\eta}^{-2}(0.5)\sum_{t=1}^{n}\left[  \Delta x_{t}\right]
\left[  x_{t-1}\right]  \overset{p}{\rightarrow}-1. \tag{A19}%
\end{equation}
When $-0.5<d-d_{0}<0$, by using $(A8)$ and $(A14)$, we have%
\begin{equation}
n^{-1}\kappa_{\eta}^{-2}(\delta)\sum_{t=1}^{n}\left[  \Delta x_{t}\right]
\left[  x_{t-1}\right]  \overset{p}{\rightarrow}-(\frac{1}{2}+\delta
)\frac{\Gamma(1+\delta)}{\Gamma\left(  1-\delta\right)  }. \tag{A20}%
\end{equation}
When $d-d_{0}=0$, by using $(A9)$ and $(A15)$, we have%
\begin{equation}
n^{-1}\kappa_{\eta}^{-2}(0)\sum_{t=1}^{n}\left[  \Delta x_{t}\right]  \left[
x_{t-1}\right]  \Rightarrow-\frac{1}{2}\left\{  \mathbf{w}^{2}(1)-1\right\}  .
\tag{A21}%
\end{equation}
When $0<d-d_{0}<0.5$, by using $(A10)$ and $(A16)$, we have%
\begin{equation}
n^{-1-2\delta}\kappa_{\eta}^{-2}(\delta)\sum_{t=1}^{n}\left[  \Delta
x_{t}\right]  \left[  x_{t-1}\right]  \Rightarrow\frac{1}{2}\mathbf{w}%
_{\delta}^{2}(1). \tag{A22}%
\end{equation}
When $d-d_{0}=m+0.5$, $m\geq1$, by using $(A11)$ and $(A17)$, we have%
\begin{equation}
n^{-2(m+1)}\left(  \log^{-1}n\right)  \kappa_{\eta}^{-2}(0.5)\sum_{t=1}^{n}
\left[  \Delta x_{t}\right]  \left[  x_{t-1}\right]  \Rightarrow\frac{1}%
{2}\mathbf{w}_{0.5,m+1}^{2}(1). \tag{A23}%
\end{equation}
When $d-d_{0}=m+\delta$, $m\geq1$%
\begin{equation}
n^{-1-2(m+\delta)}\kappa_{\eta}^{-2}(\delta)\sum_{t=1}^{n}\left[  \Delta
x_{t}\right]  \left[  x_{t-1}\right]  \Rightarrow\frac{1}{2}\mathbf{w}%
_{\delta,m+1}^{2}(1). \tag{A24}%
\end{equation}
Hence, using respectively (A1,A19), (A2,A20), (A3,A21), (A4,A22), (A5,A23),
(A6,A24) and the continuous mapping theorem, we obtain When $d-d_{0}=-0.5,$%
\begin{equation}
\left(  \log n\right)  \widehat{\rho}_{n}=\frac{n^{-1}\kappa_{\eta}%
^{-2}(0.5)\sum_{t=1}^{n}\left[  \Delta x_{t}\right]  \left[  x_{t-1}\right]
}{n^{-1}\left(  \log^{-1}n\right)  \kappa_{\eta}^{-2}(0.5)\sum_{t=1}%
^{n}\left[  x_{t-1}\right]  ^{2}}\Rightarrow\frac{-1}{\int_{0}^{1}%
\mathbf{w}_{0.5}^{2}(r)dr}. \tag{A25}%
\end{equation}
When $-0.5<d-d_{0}<0,$%
\begin{equation}
n^{1+2\delta}\widehat{\rho}_{n}=\frac{n^{-1}\kappa_{\eta}^{-2}(\delta
)\sum_{t=1}^{n}\left[  \Delta x_{t}\right]  \left[  x_{t-1}\right]
}{n^{-2-2\delta}\kappa_{\eta}^{-2}(\delta)\sum_{t=1}^{n}\left[  x_{t-1}%
\right]  ^{2}}\Rightarrow\frac{-(\frac{1}{2}+\delta)\frac{\Gamma(1+\delta
)}{\Gamma\left(  1-\delta\right)  }}{\int_{0}^{1}\mathbf{w}_{\delta}^{2}%
(r)dr}. \tag{A26}%
\end{equation}
When $d-d_{0}=0,$%
\begin{equation}
n\widehat{\rho}_{n}=\frac{n^{-1}\kappa_{\eta}^{-2}(0)\sum_{t=1}^{n}\left[
\Delta x_{t}\right]  \left[  x_{t-1}\right]  }{n^{-2}\kappa_{\eta}^{-2}%
(0)\sum_{t=1}^{n}\left[  x_{t-1}\right]  ^{2}}\Rightarrow\frac{\frac{1}%
{2}\left[  \mathbf{w}^{2}(1)-1\right]  }{\int_{0}^{1}\mathbf{w}^{2}(r)dr}.
\tag{A27}%
\end{equation}
When $0<d-d_{0}<0.5,$%
\begin{equation}
n\widehat{\rho}_{n}=\frac{n^{-1-2\delta}\kappa_{\eta}^{-2}(\delta)\sum
_{t=1}^{n}\left[  \Delta x_{t}\right]  \left[  x_{t-1}\right]  }%
{n^{-2-2\delta}\kappa_{\eta}^{-2}(\delta)\sum_{t=1}^{n}\left[  x_{t-1}\right]
^{2}}\Rightarrow\frac{\frac{1}{2}\mathbf{w}_{\delta}^{2}(1)}{\int_{0}%
^{1}\mathbf{w}_{\delta}^{2}(r)dr}. \tag{A28}%
\end{equation}
When $d-d_{0}=m+0.5,$ $m\geq1,$%
\begin{equation}
n\widehat{\rho}_{n}=\frac{n^{-2(m+1)}\left(  \log^{-1}n\right)  \kappa_{\eta
}^{-2}(0.5)\sum_{t=1}^{n}\left[  \Delta x_{t}\right]  \left[  x_{t-1}\right]
}{n^{-2\left(  m+1.5\right)  }\left(  \log^{-1}n\right)  \kappa_{\eta}%
^{-2}(0.5)\sum_{t=1}^{n}\left[  x_{t-1}\right]  ^{2}},\nonumber
\end{equation}
and then%
\[
n\widehat{\rho}_{n}\Rightarrow\frac{\frac{1}{2}\mathbf{w}_{0.5,m+1}^{2}%
(1)}{\int_{0}^{1}\mathbf{w}_{0.5,m+1}^{2}(r)dr}\equiv\rho_{1,\infty}.
\]
when $d-d_{0}=m+\delta$, $m\geq1,$%
\begin{equation}
n\widehat{\rho}_{n}=\frac{n^{-2(m+\delta)}\kappa_{\eta}^{-2}(\delta)\sum
_{t=1}^{n}\left[  \Delta x_{t}\right]  \left[  x_{t-1}\right]  }{n^{-2\left(
m+1+\delta\right)  }\kappa_{\eta}^{-2}(\delta)\sum_{t=1}^{n}\left[
x_{t-1}\right]  ^{2}}\Rightarrow\frac{\frac{1}{2}\mathbf{w}_{\delta,m+1}%
^{2}(1)}{\int_{0}^{1}\mathbf{w}_{\delta,m+1}^{2}(r)dr}\equiv\rho_{2,\infty}.
\tag{A30}%
\end{equation}
Now consider the $t$-statistic. First notice that\newline$\widehat{\sigma}%
_{n}^{2}=n^{-1}\left(  \sum_{t=1}^{n}\left(  \Delta x_{t}\right)
^{2}+\widehat{\rho}_{n}^{2}\sum_{t=1}^{n}\left(  x_{t-1}\right)
^{2}-2\widehat{\rho}_{n}\sum_{t=1}^{n}\left[  \Delta x_{t}\right]  \left[
x_{t-1}\right]  \right)  $. Hence, when $d-d_{0}=-0.5$, by using $(A1)$,
$(A13) $ $(A19)$ and $(A25)$, it follows%
\begin{equation}
\widehat{\sigma}_{n}^{2}\overset{p}{\rightarrow}var(\Delta x_{t}%
)=\frac{4\sigma_{\varepsilon}^{2}}{\pi}. \tag{A31}%
\end{equation}
When $-0.5<d-d_{0}<0$, by using $A2$, $A14$, $A20$ and $A26$, it follows%
\begin{equation}
\widehat{\sigma}_{n}^{2}\overset{p}{\rightarrow}var(\Delta x_{t})=\frac
{\sigma_{\varepsilon}^{2}\Gamma(1-2\delta)}{\Gamma^{2}(1-\delta)}. \tag{A32}%
\end{equation}
When $d-d_{0}=0$, by using $A3$, $A15$, $A21$ and $A27$, it follows%
\begin{equation}
\widehat{\sigma}_{n}^{2}\overset{p}{\rightarrow}var(\Delta x_{t}%
)=\sigma_{\varepsilon}^{2}. \tag{A33}%
\end{equation}
When $0<d-d_{0}<0.5$, by using $A4$, $A16$, $A22$ and $A28$, it follows%
\begin{equation}
\widehat{\sigma}_{n}^{2}\overset{p}{\rightarrow}var(\Delta x_{t})=\frac
{\sigma_{\varepsilon}^{2}\Gamma(1-2\delta)}{\Gamma^{2}(1-\delta)}. \tag{A34}%
\end{equation}
When $d-d_{0}=m+0.5$, $m\geq1$, by using $A5$, $A17$, $A23$ and $A29$, it
follows
\begin{equation}
\frac{\kappa_{\eta}^{-2}(0.5)\widehat{\sigma}_{n}^{2}}{n^{2m}\left(  \log
n\right)  }\Rightarrow\int_{0}^{1}\mathbf{w}_{0.5,m}^{2}(r)dr+\rho_{1,\infty
}^{2}\int_{0}^{1}\mathbf{w}_{0.5,m+1}^{2}(r)dr-\rho_{1,\infty}\mathbf{w}%
_{0.5,m+1}^{2}(1). \tag{A35}%
\end{equation}
When $d-d_{0}=m+\delta$, $m\geq1$ by using $A6$, $A18$, $A24$ and $A30$, it
follows
\begin{equation}
n^{-2m-2\delta+1}\kappa_{\eta}^{-2}(\delta)\widehat{\sigma}_{n}^{2}%
\Rightarrow\int_{0}^{1}\mathbf{w}_{\delta,m}^{2}(r)dr+\rho_{2,\infty}^{2}%
\int_{0}^{1}\mathbf{w}_{\delta,m+1}^{2}(r)dr-\rho_{2,\infty}\mathbf{w}%
_{\delta,m+1}^{2}(1). \tag{A36}%
\end{equation}
Finally, by using respectively (A1,A19,A31), (A2,A20,A32), (A3,A21,A33),
(A4,A22,A34) (A5,A23,A35), (A6,A24,A36) we obtain for the $t$%
-statistic:\newline\noindent When $d-d_{0}=-0.5,$ $t_{\widehat{\rho}_{n}%
}=\frac{n^{-0.5}(\log^{-0.5}n)\kappa_{\eta}^{-1}(0.5)\sum_{t=1}^{n}\left[
\Delta x_{t}\right]  \left[  x_{t-1}\right]  }{\widehat{\sigma}_{n}\left\{
n^{-1}(\log^{-1}n)\kappa_{\eta}^{-2}(0.5)\sum_{t=1}^{n}\left[  x_{t-1}\right]
^{2}\right\}  ^{0.5}}\overset{p}{\rightarrow}-\infty,$ and $t_{\widehat{\rho
}_{n}}=O_{p}(n^{0.5}\log^{0.5}n)$.\newline\noindent When $-0.5<d-d_{0}<0,$
$t_{\widehat{\rho}_{n}}=\frac{n^{-\delta-1}\kappa_{\eta}^{-1}(\delta
)\sum_{t=1}^{n}\left[  \Delta x_{t}\right]  \left[  x_{t-1}\right]  }%
{\widehat{\sigma}_{n}\left\{  n^{-2\delta-2}\kappa_{\eta}^{-2}(\delta
)\sum_{t=1}^{n}\left[  x_{t-1}\right]  ^{2}\right\}  ^{0.5}}\overset
{p}{\rightarrow}-\infty,$ and $t_{\widehat{\rho}_{n}}=O_{p}(n^{-\delta}%
).$\newline\noindent When $d-d_{0}=0,$ $t_{\widehat{\rho}_{n}}=\frac
{n^{-1}\kappa_{\eta}^{-1}(0)\sum_{t=1}^{n}\left[  \Delta x_{t}\right]  \left[
x_{t-1}\right]  }{\widehat{\sigma}_{n}\left\{  n^{-2}\kappa_{\eta}^{-2}%
(0)\sum_{t=1}^{n}\left[  x_{t-1}\right]  ^{2}\right\}  ^{0.5}}\Rightarrow
\frac{\frac{1}{2}\left[  \mathbf{w}^{2}(1)-1\right]  }{\left[  \int_{0}%
^{1}\mathbf{w}^{2}(r)dr\right]  ^{0.5}},$ and $t_{\widehat{\rho}_{n}}%
=O_{p}(1).$\newline\noindent When $0<d-d_{0}<0.5,$ $t_{\widehat{\rho}_{n}%
}=\frac{n^{-\delta-1}\kappa_{\eta}^{-1}(\delta)\sum_{t=1}^{n}\left[  \Delta
x_{t}\right]  \left[  x_{t-1}\right]  }{\widehat{\sigma}_{n}\left\{
n^{-2\delta-2}\kappa_{\eta}^{-2}(\delta)\sum_{t=1}^{n}\left[  x_{t-1}\right]
^{2}\right\}  ^{0.5}}\overset{p}{\rightarrow}-\infty,$ and $t_{\widehat{\rho
}_{n}}=O_{p}(n^{\delta}).$\newline\noindent When $d-d_{0}=m+0.5,$ $m\geq1,$
\newline$t_{\widehat{\rho}_{n}}=\frac{n^{0.5}\left\{  n^{-2(m+1)}\left(
\log^{-1}n\right)  \kappa_{\eta}^{-2}(0.5)\sum_{t=1}^{n}\left[  \Delta
x_{t}\right]  \left[  x_{t-1}\right]  \right\}  }{\left\{  n^{-m}\left(
\log^{-0.5}n\right)  \kappa_{\eta}^{-1}(0.5)\widehat{\sigma}_{n}\right\}
\left\{  n^{-2\left(  m+1.5\right)  }\left(  \log^{-1}n\right)  \kappa_{\eta
}^{-2}(0.5)\sum_{t=1}^{n}\left[  x_{t-1}\right]  ^{2}\right\}  ^{0.5}}%
\overset{p}{\rightarrow}+\infty,$ and $t_{\widehat{\rho}_{n}}=O_{p}(n^{0.5}%
).$\newline\noindent When $d-d_{0}=m+\delta$, $m\geq1,$\newline$t_{\widehat
{\rho}_{n}}=\frac{n^{0.5}\left\{  n^{-1-2(m+1)}\kappa_{\eta}^{-2}(\delta
)\sum_{t=1}^{n}\left[  \Delta x_{t}\right]  \left[  x_{t-1}\right]  \right\}
}{\left\{  n^{-m-\delta+0.5}\kappa_{\eta}^{-1}(\delta)\widehat{\sigma}%
_{n}\right\}  \left\{  n^{-2\left(  m+1+\delta\right)  }\kappa_{\eta}%
^{-2}(\delta)\sum_{t=1}^{n}\left[  x_{t-1}\right]  ^{2}\right\}  ^{0.5}%
}\overset{p}{\rightarrow}+\infty,$ and $t_{\widehat{\rho}_{n}}=O_{p}%
(n^{0.5}).$

\bigskip

\textbf{Appendix 2: }$\mathbf{F-DF}$\textbf{\ test in the Presence of
deterministic components}

We assume that the univariate process $\left\{  y_{t},\text{ }t\in%
\mathbb{Z}
\right\}  $ can be generated by the following two mechanisms%
\begin{equation}
y_{t}=\alpha+(1-L)^{-d}u_{t},\text{ \ \ }t\in%
\mathbb{Z}
, \tag{A}%
\end{equation}
and%
\begin{equation}
y_{t}=\alpha+\beta t+(1-L)^{-d}u_{t},\text{ \ \ }t\in%
\mathbb{Z}
, \tag{B}%
\end{equation}
where $u_{t}$ is as in $(2.1)$. For $DGP$ $\left(  A\right)  $ and $(B)$ we
can use the $F-DF$ test without having to use their specific asymptotic
theory. This can be done by differencing the process $y_{t}$ one time in $(A)
$ and twice in $(B)$. To be more clear, we consider testing hypotheses%
\begin{equation}
H_{0}:d\geq d_{0}\text{ \ \ against \ \ }H_{1}:d<d_{0}, \tag{C}%
\end{equation}
even though, the framework provided in this paper does not allow us to use the
$DGP$ $(A)$ and $(B)$, we can , nevertheless use it as the following.

For the $DGP$ $(A)$, the constant, $(\alpha)$, can be removed by first
differencing,%
\[
(1-L)y_{t}=(1-L)^{1-d}u_{t},
\]
and then, by transforming the $DGP$ $(A)$ in this way, instead $(C)$, we
consider testing hypotheses%
\[
H_{0}:d_{1}\geq d_{0}-1\text{ \ \ }against\text{ \ \ }H_{1}:d_{1}%
<d_{0}-1,\text{ \ \ \ where }d_{1}=d-1,
\]
by using the auxiliary regression model%
\[
\Delta^{d_{0}-1}y_{t}=\rho\Delta^{d_{0}}y_{t-1}+\varepsilon_{t}.
\]

For the $DGP$ $(B)$, the constant, $\left(  \alpha\right)  $, and the
parameter time trend, $\left(  \beta\right)  $, can be removed by differencing
twice the process $y_{t}$,%
\[
(1-L)^{2}y_{t}=(1-L)^{d-2}u_{t},
\]
and then, for this transformed model, instead $(C)$, we consider testing
hypotheses%
\[
H_{0}:d_{2}\geq d_{0}-2\text{ \ \ }against\text{ \ \ }H_{1}:d_{2}%
<d_{0}-2,\text{ \ \ \ where }d_{2}=d-2,
\]
by using the auxiliary regression model%
\[
\Delta^{d_{0}-2}y_{t}=\rho\Delta^{d_{0}-1}y_{t-1}+\varepsilon_{t}.
\]

\end{document}